\documentclass[12pt]{UOthesis}

\usepackage{fancybox}
\usepackage{shadow}

\allowdisplaybreaks

\NoLogo
\NoSignatures
\NoResume
\NoDedications
\NoListoftables
\NoListoffigures

\setUOname{Verne Cazaubon}         
\setUOcpryear{2008}               
\setUOtitle{In search of a Lebesgue density theorem for $\mathbb{R}^{\infty}$}  

\msc               

\setUOabstract{

We look at a measure, $\lambda^\infty$, on the infinite-dimensional space, $\RR^\infty$, for which we attempt to put forth an analogue of the Lebesgue density theorem.  Although this measure allows us to find partial results, for example for continuous functions, we prove that it is impossible to give an analogous theorem in full generality.  In particular, we proved that the Lebesgue density of probability density functions on $\RR^\infty$ is zero almost everywhere.
}

\setUOresume{
  Vous pouvez introduire le r\'esum\'e de votre th\`ese ici.
}

\setUOthanks{
There are many people that I would like to thank for helping me get through these two years from September 2006 to October 2008.

First and foremost, I thank God.  Without Him I would not have had the patience or strength to get through this work.  He has guided me at every step and made sure that the burden was never more than I could handle.

I am tremendously grateful to the Canadian Government and in particular the Canadian Bureau for International Education (CBIE) for choosing me to receive a two year full scholarship to attend an institute of my choice in Canada.  Without that I would not have been here.  I appreciate all that you have done and I am forever indebted to you for what you have given me is an education, and that is priceless.
	
I would like to thank my supervisor, Vladimir Pestov, for his guidance throughout these trying times, for always keeping me on track since it was so easy to go off on a tangent given the enormous scope of this topic, and for his contagious enthusiasm when it comes down to talking Mathematics.  For his patience with me since I had to learn so much to get to a fraction of his level.  For his interest in other things apart from this	thesis and for his advice on non-school related matters.

I thank my examiners Wojciech Jaworski and David McDonald for their input and for making sure that this work is up to standard with others in the Mathematics community.  I am more confident now that it has been scrutinized by some of the best minds in the area and met their approval.
	
A big thank you to my family back home in St. Lucia for their support and prayers.  They have always encouraged me to strive for the best and without their words of wisdom, I may not have been at this level today.  I have not seen them in over two years but I know that even though they are physically over two thousand miles away, they are right here with me.
	
I thank my brother, Luke Cazaubon, for his encouraging words and for allowing me to stay at his place for a year and a half, 	providing shelter, food, needs and wants.  He always asked that all I do is focus on getting this thesis done, and worry about nothing else, so this final work is dedicated to him also.

Also, I owe a great deal of gratitude to Janet Maria Prins, who was and is always encouraging me to be the best that I can be.  For the kindness and the hospitality of her family and for those weekend getaways to Pembroke which were always very relaxing, thanks.

I would also like to thank the office mates of room 311 at the Mathematics Department at uOttawa, whose graduate students 		always welcomed me in there with open arms.  Thanks to Anjayan, Rachelle, Trevor and Vladan who always had some interesting 		discussion going on.  I would especially like to thank Wadii Hajji for the many discussions we had on the topic, and for his 	insight and his questioning which allowed me to see things in a different (and sometimes clearer) light. 

And to Stephen Manners, whose generosity kept me comfortable during the last eight months of my stay in Canada; thank you for helping me out when I was in a tight situation.  To all who I forgot to mention and everyone else who in the smallest way contributed to the completion of this thesis, I am very appreciative. \\

\begin{flushright}
Verne Cazaubon\\
Ottawa, Canada, October 2008
\end{flushright}
}

\setUOdedicationsText{
  In memory of my famous bottle of port which was taken away from me
  just after reaching full maturity.
}

{\theorembodyfont{\rm}
  \newtheorem{rmk}[theo]{Remark}
  \newtheorem{egg}[theo]{Example}
}

\newcommand{\RR}{\mathbb{R}}
\newcommand{\NN}{\mathbb{N}}
\newcommand{\II}{\mathbb{I}}

\DeclareMathOperator{\esssup}{\text{ess sup}}

\begin{document}


\chapter{Introduction}

The aim of the thesis at hand is to analyse the possible extension of the Lebesgue density theorem to infinite dimensions.  We will attempt to explain why it is natural and interesting to ask this question.  The direction in which we were planning to take the original research project is quite different from the results which we now present and here we will explain where the motivation comes from.

\section{Motivation}

\subsection{Similarity Search}

The idea which grew into this thesis was the following: investigate the existence of an efficient algorithm which, given a point $q$ known as a query point, will search a high-dimensional object (call it $D$ for database) and return the nearest point or the nearest points to $q$.  This is the basic idea behind the well known similarity search, also known as neighbour search and closest point search, and in the case where we accept more than one output point it is also known as $k$-nearest neighbour search and proximity search.  For a good introduction to high-dimensional similarity search, see chapter 9 in \cite{SS}.  Similarity search has applications in coding theory, database and data mining, statistics and data analysis, information retrieval, machine learning and pattern recognition amongst other fields.   If these applications mean nothing to you, some more specific applications would be searching a database of resum\'es for an applicant who suits a particular job description, searching for a used car by price, history, mileage, make, model, year and colour and even searching for a companion through an online dating service.

Each of our points (entries) in the database is described by a number of features, and this number may range from a small number like two or three, which means that our objects are described with little detail, to hundreds of thousands and even millions, which means our objects are extremely detailed.  For many applications, we may not be able to get a match satisfying all our desired features, in which case it would suffice to return matches satisfying most of the features.  This may be decided by giving different weights to the features and choosing matches with the highest weight.  This is not necessary if all features have the same priority.  Also, we may not want only one match.  We may want our algorithm to return a number of matches and leave it up to the human to choose the best match according to their own discretion, personal preference, or according to some criteria which is not taken into consideration in the database.  For example, when searching for a companion, what one person finds visually appealing is personal and there is no criteria which will allow an algorithm to pick out an exact match.  In these cases, $k$-nearest neighbour search is preferred over nearest neighbour search.

Now, a very relevant question would be, what do we mean by a high-dimensional object?  Well, this object is usually a space in which our data points lie.  The dimension, $d$, of our space depends on the number of features of each point.  It is greater than or equal to the number of features of the point with the highest number of features.  For example, if our data consists of information about an employee at some company, then an entry may include features such as the worker's employee number, first name, last name, gender, position at the company and maybe even marital status.  Thus, a 6-dimensional database is enough to store information about employees.  So, the question can now be stated as follows: given a database of employees, $D$, with the features mentioned above stored in $D$, and some data, $q$, about someone, find the employee in $D$ who is the closest match to $q$, or return the $k$ workers who best fit the description.

For simplicity, and since we are studying the asymptotic case, we choose as our domain $\RR^\infty$ to be the set of all possible points.  A typical example of a domain would be a Euclidean space with Euclidean distance or some $\ell_p$ metric.  There are many algorithms to search through data structures like k-d trees, M-trees and locality-sensitive hashing (again see \cite{SS} for more examples and for references to literature as this work has veered away from the topic).  A model that we wanted to investigate uses a random geometric graph and we assumed a greedy walk type algorithm.  Start at a random vertex.  Take one step to another node which is closer to the query.  From this new vertex, repeat the last step and continue until we can get no closer to the query point, then stop.  We do not deal with the dynamic case where nodes are removed, instead, in our case, new nodes are added but the older ones remain.  The most difficult concept to grasp is the assumption that as our dimension grows, and more is known about each object, our overview of them changes only slightly in that as we see more data we expect it to agree in some sense to what we have already seen.

In the case where data lies in a low-dimensional space, very efficient algorithms exist to find the solution to our query \cite{HE}.  However, as dimension increases 
algorithms fall prey to what is known as the ``curse of dimensionality'' (see \cite{SS} for a discussion).  Informally, this means that the higher our dimension, the more difficult it is to find the nearest neighbour(s) of a given point.  In fact, given a high enough dimension, known algorithms perform no better than a brute force search \cite{WSB}.  It is widely believed that in high-dimensional databases there are no efficient search algorithms, but it remains unproven at a mathematical level, so we are interested in a framework for the analysis of this problem.

In an effort to solve this problem of inefficiency in high dimensions, it would be nice to find an algorithm which is independent of dimension.  Of course, this is impossible, so we settle for something a bit more realistic, and that is an algorithm with time complexity of order $ O(\log d)$.  Now, the amount of work which goes into searching our database depends on how our data is stored.  If data is stored haphazardly, then the time to find the nearest neighbour(s) of a query point may increase exponentially.  Since we are concerned with efficiency, we have decided that the most suitable structure to model our data would be a random graph, but not just any random graph, a random geometric graph.

\subsection{Random Geometric Graphs}

Whenever a mathematician hears the words `random graph', the first model which comes to mind is the Erd\H os-R\'enyi model, $G(n,p)$, where we start with $n$ vertices and choose each of the edges of the graph with equal probability $p$.  Whilst the study of these kinds of graphs is useful, they are not very realistic and in this situation are of little use to us.  If we attempt to search through this model of a random graph, then there is no assurance that `jumping' from one vertex to another will bring us any closer to our target, $q$.  However, random graphs of this type are still useful for proving the existence of graphs satisfying certain properties.

In our investigation, we planned to study a different kind of random graph model which has been of interest recently.  These graphs have a little more structure than those of the Erd\H os-R\'enyi model, and it is this structure which we will try to take advantage of.  By their very nature, they are an excellent choice to model data in a high-dimensional space.  They are still random, and no less random than any other model, but random in a sense which is more suitable to our situation.  Random geometric graphs are formed by arbitrarily choosing $n$ points in a space with regard to some probability distribution, and joining any two points, separated by less than some specified distance, with an edge.  This model is used to study distributed wireless networks, sensory based communication networks, Percolation Theory and cluster analysis --- which is a powerful technique with applications in Medicine, Biology and Ecology.  An excellent reference which introduces and goes in-depth into random geometric graphs is \cite{MP}.

Our initial aim was to use random geometric graphs in the study of similarity search.  Of course, random geometric graphs is an obvious model on which such work can be done.  If given a set (database) of $n$ points in a $d$-dimensional metric space and a specific query point $q$, we would like to find the nearest neighbour, $p$, to that point $q$.  When using random geometric graphs, we would connect all points which are within a distance $r$ from each other, then choose any point to start our search.  By using the greedy algorithm described earlier, we would get a result.  The asymptotic setting in which we would like to work and in which researchers model similarity search is by taking both $n$ and $d$ tending toward infinity.  This setting requires us to work in $\RR^{\infty}$.

\section{Our Setting}

Let $f$ be a probability density function on $\RR^{\infty}$ which is bounded and integrable with regard to a certain measure, $\lambda^{\infty}$, that is, $\int_{\RR^{\infty}} f \, d \lambda^{\infty} = 1$.  Let us assume that our data is modelled by a sequence $X_1, X_2, \ldots$ which are independent and identically distributed random variables taking their values in $\RR^\infty$.  
Let $\RR^{\infty}$ be equipped with a metric, $\rho(\cdot, \cdot)$, which induces the product topology.
\begin{egg}
	The metric
	\[
		\rho(x,y) = \sum_{i=1}^{\infty} 2^{-i} \frac{|x_i - y_i|}{1 + |x_i - y_i|}
	\]
	where $x = (x_i), y = (y_i) \in \RR^{\infty}$, induces the product topology on $\RR^{\infty}$.
\end{egg}
Let $\mathfrak{X}_n = \{ X_1, X_2, \ldots, X_n \} \subset \RR^{\infty}$ be the vertex set of an undirected graph $G(\mathfrak{X}_n, r)$ where any two vertices $X_i, X_j$ are joined if and only if $\rho(X_i, X_j) < r$, where $r$ denotes the radius.  We call $G(\mathfrak{X}_n, r)$ a random geometric graph.

Now if we let our random geometric graph have a fixed radius $r$, but continue to add vertices (increase $n$ while keeping $d$ constant), then the average degree of each vertex is guaranteed to rise.  This is undesirable if we are using a greedy algorithm to search our graph because it increases search time.  To curb this phenomenon and others similar to it, we introduce a sequence of radii $(r_n)$ and limiting regimes \cite{MP}.  We are interested in one particular limiting regime for $(r_n)$, that is the \emph{thermodynamic limit}.  In this limiting regime, the expected degree of a typical vertex tends to be constant.  Here, $r_n \propto cn^{-\frac{1}{d}}$, for some constant $c$.  It can be shown that if the limiting constant is taken to be above some critical value, then with high probability a giant component will arise in $G(\mathfrak{X}_n; r_n)$.  This giant component is a connected subgraph which contains most of the vertices of the entire graph.  It is studied in \cite{CL, MP}.  There are other limiting regimes which we may come across, these include:
\begin{itemize}
	\item[(i)]		the sparse limit regime: $nr_n^d \rightarrow 0$
	\item[(ii)]		the dense limit regime: $nr_n^d \rightarrow \infty$
	\item[(iii)]	the connectivity regime: $r_n \propto c((\log n)/n)^{-\frac{1}{d}}$ (this is a special case of the dense limit 										regime).
\end{itemize}

Knowing these limits is of extreme importance, and here is why.  If we are looking for the nearest neighbours to our query point $q$, we would hope that $q$ is not part of the giant component, because our search becomes much easier when we have eliminated a large portion of vertices.  Another scenario which may be of importance is finding a subgraph in our model which is isomorphic to some other graph $\Gamma$.  We will look at the latter in more detail as it will propel us to our main objective.

We have found that proving analogues of certain results about the behaviour of random geometric graphs in $\RR^d$ as $d$ tends to infinity requires using the Lebesgue density theorem in $\RR^\infty$, and so we have concentrated on the task of finding an analogue of this result.

Let us say a word about the so called Lebesgue measure, $\lambda^{\infty}$ on $\RR^\infty$.  It is known that given any infinite-dimensional locally convex topological vector space, $X$, there does not exist a non-trivial translation-invariant sigma-finite Borel measure on $X$ \cite{RB, HSY}.  In other words, Haar measure does not exist in an infinite-dimensional setting.  However, there have been a few mathematicians who have written of analogues to the Lebesgue Measure in infinite dimensions, they include R. Baker \cite{RB}, Y. Yamasaki \cite{YY}, and the trio of N. Tsilevich, A. Vershik and M. Yor \cite{TVY}.  In \cite{RB}, a nontrivial translation invariant Borel measure, $\lambda^{\infty}$, on $\RR^{\infty}$ is introduced.  This is a sigma-additive Borel measure which is analogous to the Lebesgue measure, but it is of course not sigma-finite.  We focus on this measure in the hope that it will lead us to a positive result.

An interesting circumstance of both $n$ and $d$ tending to infinity is that at any point in time, only the first $d$ coordinates (features) are revealed to us while the rest are kept hidden.  As time goes by, however, in addition to obtaining new data points, our knowledge of the data increases as we discover new features and their values at previous data points.  Thus, at any step in time we have $n$ points modelled by independent and identically distributed random variables distributed with respect to $f_d \cdot \lambda^d$, where $f_d$ and $\lambda^d$ can be thought of as projection of our probability density function, $f$, and the push forward of our measure, $\lambda^{\infty}$, to the finite-dimensional case.  Now, these functions, $f_d$, depend only on the first $d$ coordinates and as we get to know more information, our functions do not vary greatly.  This is intuitively what we want.  We do not want or expect the data that is obtained in the future to be extremely different from what we have at any point in time.  
Thus, we are not dealing with complete independence.  And let us not forget that $n$ is also increasing, but $n$ is tied to $d$ in the sense that it is not independent.  It grows faster than $d$, but only sub-exponentially.


More precise definitions of $f_d$ and $\lambda^{\infty}$ will be given later on.  $\lambda^d$ is just the $d$-dimensional Lebesgue measure which is well known.  What we should keep in mind is that at any instant in time, we have a random geometric graph in $\RR^d$.

\section{Overview}



The plan for the rest of the paper is as follows: In chapter 2, we provide the base and background needed to continue through the paper without getting lost.  The terminology and notation introduced there will be used throughout the rest of the paper.  New or specific terms will be defined in the chapter in which they are put forth.

Chapter 3 explains the reason for embarking on this project.  We go through the theorem which started it all.  We see for the first time the Lebesgue density theorem, the extension of which is the purpose of this paper.  We review briefly the concept of a Lebesgue point and state Lebesgue's density theorem.  We see why it is not possible to simply extend this concept, as is, to an infinite-dimensional case.

In chapter 4, we deal with the measure which we will be using in our analysis.  It is a very interesting measure and it is supposedly very nice to work with seeing it is analogous to the Lebesgue measure.  We start by showing that a non-trivial translation-invariant locally finite measure is not possible in an infinite-dimensional setting.  This is a weaker statement than that proved by Andr\'e Weil, but still shows us that working in infinite dimensions is not easy.  We then go through the construction of the ``Lebesgue measure'' on $\RR^\infty$.

Chapter 5 is concerned with a not very well known example by Jean Dieudonn\'e.  It was instrumental in shaping the outcome of this paper.  We go through his paper in detail as this is probably the first English translation of the result, after which we point out the complication which arose due to the result.

Chapter 6 introduces a new class of functions which satisfy the analogous Lebesgue density theorem and brings to the fore a very interesting theorem by Jessen.  This theorem was crucial in attempting to solve our problem.  It is just as interesting, and complex, that the generalisation of this theorem does not hold as shown by Dieudonn\'e in chapter 5.

The final chapter exhibits some routes on which we embarked but departed from for some reason or the other.  It also shows some of the difficulties we had in general.  In the last section we produce the final result showing the impossibility of the density theorem on $\RR^\infty$, in full generality as required by our original goals.

The last detail we would like to point out before continuing is that double daggers $(\ddag)$ will indicate new results.




 



\cleardoublepage


\chapter{Preliminaries}

\section{Some Important Definitions}

The purpose of this chapter is to introduce some of the terms (and general notation) that will be used throughout the paper.  As we come across or introduce specific terms, they will be defined in the chapter in which they appear.  The definitions here have been compiled from \cite{RD, HP, RHK, JLK, SL, JEM, JRM, CAR, WR}.

\subsection{Topology}

A collection $\tau$ of subsets of a set $X$ is said to be a \emph{topology} in $X$ if $\tau$ has the following three 						properties:
\begin{itemize}
	\item $\emptyset \in \tau$ and $X \in \tau$
	\item If $A_i \in \tau$ for $i = 1, \ldots, n$, then $A_1 \cap A_2 \cap \cdots \cap A_n \in \tau$
	\item If $\{A_{\alpha}\}$ is an arbitrary collection of members of $\tau$, then $\cup_{\alpha} A_{\alpha} \in \tau$.
\end{itemize}

\begin{egg}
	Let $X_1 = \{1, 2, 3, 4\}$, then the following are topologies of $X_1$: \\
	$\tau_1 = \{\emptyset, \{1\}, \{2\}, \{1,2\}, \{1,2,3,4\}\}$ \\
	$\tau_2 = \{\emptyset, \{3\}, \{1,2,3\}, \{1,2,3,4\}\}$
\end{egg}

If $\tau$ is a topology on $X$, then $(X, \tau)$ is a \emph{topological space}, and the members of $\tau$ are called \emph{open sets} in $X$.  A set $B \subseteq X$ is \emph{closed} if its complement $B^c$ is open.  Both $\emptyset$ and $X$ are closed, finite unions of closed sets are closed and arbitrary intersections of closed sets are closed.  A set can be both open and closed and a set can be neither open nor closed.  $X$ is called \emph{connected} if and only if the only subsets of $X$ which are both open and closed in $\tau$ are the empty set, $\emptyset$, and $X$ itself.  The \emph{closure} of a set $B \subseteq X$ is the smallest closed set in $X$ which contains $B$.

\begin{egg}
	The sets $\emptyset$ and $X_1$ are always both open and closed.  $\{3\} \subset X_1$ is neither open nor closed with respect 		to $\tau_1$, however, it is an open set with respect to $\tau_2$.  The closure of $\{3\}$ with respect to $\tau_2$ is $X_1$.
\end{egg}

A set $U$ in a topological space $(X, \tau)$ is a \emph{neighbourhood} of a point $x$ if and only if $U$ contains an open set to which $x$ belongs.  A neighbourhood of a point need not be an open set, but every open set is a neighbourhood of each of its points.  A family $\mathcal{B}$ of sets is a \emph{base} for a topology $\tau$ if and only if $\mathcal{B}$ is a subfamily of $\tau$ and for each point $x$ of the space, and each neighbourhood $U$ of $x$, there is a member $V$ of $\mathcal{B}$ such that $x \in V \subset U$.  A space whose topology has a countable base is called \emph{second countable}.  Any second countable space is \emph{separable} but not vice versa.  A family $\mathcal{S}$ of sets is a \emph{subbase} for a topology $\tau$ if and only if the class of finite intersections of members of $\mathcal{S}$ is a base for $\tau$.  Equivalently, iff each member of $\tau$ is the union of finite intersections of the members of $\mathcal{S}$, then $\mathcal{S}$ is a subbase for $\tau$.  A map $f$ of a topological space $(X, \tau)$ into a topological space $(Y, \upsilon)$ is \emph{continuous} if and only if $f^{-1}(U) \in \tau$ for each $U \in \upsilon$.  That is, iff the inverse of each open set is open then $f$ is continuous.

A family $\mathcal{C}$ of sets is a \emph{cover} of a set $X$ if and only if $X$ is a subset of the union $\bigcup\{A:A \in \mathcal{C}\}$.  That is, iff each member of $X$ belongs to some $A \in \mathcal{C}$.  It is an open cover if and only if each $A \in \mathcal{C}$ is an open set.  A \emph{subcover} of $\mathcal{C}$ is a subset of $\mathcal{C}$ that still covers $X$.  A set $B \subset X$ is \emph{compact} if every open cover of $B$ contains a finite subcover.  If $X$ is compact, it is called a compact space.

\begin{defn} \label{D:Ideal}
Given a set $X$, a non-empty subset, $I$, of the power set of $X$ is called an \emph{ideal} on $X$ if:
\begin{itemize}
	\item $A \in I$ and $B \subseteq A$ implies $B \in I$ and
	\item $A, B \in I$ implies $A \cup B \in I$.
\end{itemize}
\end{defn}

\subsection{Metric Spaces}

The most familiar topological spaces are metric spaces.  A \emph{metric space} is a set $X$ in which a distance function (or metric) $\rho$ is defined, with the following properties:
\begin{itemize}
	\item $0 \leq \rho(x,y) < \infty$ for all $x, y \in X$
	\item $\rho(x,y) = 0$ iff $x=y$
	\item $\rho(x,y) = \rho(y,x)$ for all $x, y \in X$
	\item $\rho(x,y) \leq \rho(x,z) + \rho(z,y)$ for all $x, y, z \in X$
\end{itemize}

\begin{defn} \label{D:PositivelySeparated}
	Two sets $A, B$ in a metric space $X$ with metric $\rho$ are \emph{positively separated} if
	\[
		\inf_{\substack{a \in A\\b \in B}} \rho(a, b) > 0.
	\]
\end{defn}

Let $(X, \rho_1)$ and $(Y, \rho_2)$ be two metric spaces.  A function $f:X \rightarrow Y$ is called \emph{uniformly continuous} if for every $\epsilon > 0$, there exists a $\delta > 0$ such that for all $x_1, x_2 \in X$ with $\rho_1(x_1, x_2) < \delta$ we have $\rho_2(f(x_1), f(x_2)) < \epsilon$.  Every uniformly continuous function is continuous, but not vice versa.

\begin{defn}
Let $x_n$ be a sequence of numbers.  We say $x_n$ \emph{converges} to $x$ if for any number $\epsilon > 0$ there is an integer $N$ such that $|x_n-x| < \epsilon$ for all integers $n \geq N$.  This is written as $\lim_{n \rightarrow \infty} x_n = x$ or $x_n \rightarrow x$ as $n \rightarrow \infty$.  A sequence $x_n$ in $\RR$ is called a \emph{Cauchy sequence} if for every number $\epsilon > 0$ there is an integer $N$ (depending on $\epsilon$), such that $|x_n-x_m| < \epsilon$ whenever $n \geq N$ and $m \geq N$.
\end{defn}

\subsection{Measure Theory}

A collection $\mathfrak{M}$ of subsets of a set $X$ is said to be a \emph{$\sigma$-algebra} in $X$ if it has the following 			properties:
\begin{itemize}
	\item $X \in \mathfrak{M}$
	\item If $A \in \mathfrak{M}$, then $A^c \in \mathfrak{M}$, where $A^c$ is the complement of $A$ relative to $X$.
	\item If $A = \bigcup_{n=1}^\infty A_n$ and if $A_n \in \mathfrak{M}$ for $n = 1, 2, \ldots$, then $A \in \mathfrak{M}$.
\end{itemize}

\begin{egg}
	Let $X_1 = \{1, 2, 3, 4\}$, then the following are both $\sigma$-algebras of $X_1$: \\
	$\mathfrak{M}_1 = \{\emptyset, \{1\}, \{2,3,4\}, \{1,2,3,4\}\}$ \\
	$\mathfrak{M}_2 = \{\emptyset, \{1\}, \{2\}, \{1,2\}, \{3,4\}, \{2,3,4\}, \{1,3,4\}, \{1,2,3,4\}\}$
\end{egg}

If $\mathfrak{M}$ is a $\sigma$-algebra in $X$, then $(X, \mathfrak{M})$ is called a \emph{measurable space}, and the members of $\mathfrak{M}$ are called \emph{measurable sets} in $X$.  If $X$ is a measurable space, $Y$ is a topological space, and $f$ is a mapping of $X$ into $Y$, then $f$ is said to be \emph{measurable} provided that $f^{-1}(V)$ is a measurable set in $X$ for every open set $V$ in $Y$.

\begin{defn} \label{D:Premeasure}
	A function $\xi$ on a family $\mathfrak{C}$ of subsets of $X$ is called a pre-measure if:
	\begin{itemize}
 		\item $\emptyset \in \mathfrak{C}$ and $\xi(\emptyset) = 0$;
 		\item $0 \leq \xi(C) \leq +\infty$ for all $C$ in $\mathfrak{C}$.
	\end{itemize}
\end{defn}

\begin{defn} \label{OuterMeasure}
A function $\mu^\star$ defined on the sets of a space $X$ is called an \emph{outer measure} on $X$ if it satisfies the conditions:
	\begin{itemize}
		\item $\mu^\star(E)$ takes values in $[0, + \infty]$ for each subset $E$ of $X$;
		\item	$\mu^\star(\emptyset) = 0$;
		\item	if $E_1 \subset E_2$ then $\mu^\star(E_1) \leq \mu^\star(E_2)$; and
		\item	if $\{E_i\}$ is any sequence of subsets of $X$ then
		\[
			\mu^\star \left ( \bigcup_{i=1}^\infty E_i \right ) \leq \sum_{i=1}^\infty \mu^\star(E_i).
		\]
	\end{itemize}
\end{defn}

\begin{defn} \label{D:MetricOuterMeasure}
	An outer measure $\mu^\star$ defined on a metric space $X$ is called a \emph{metric outer measure} if
	\[
		\mu^\star(A \cup B) = \mu^\star(A) + \mu^\star(B)
	\]
	for all pairs of positively separated sets $A, B$.
\end{defn}

\begin{defn}
A \emph{measure}, $\mu$, is an outer measure restricted to a $\sigma$-algebra $\mathfrak{M}$, and which is countably additive.  This means if $\{A_i\}$ is a disjoint countable collection of members of $\mathfrak{M}$, then
\[
	\mu \left ( \bigcup_{i=1}^{\infty} A_i \right ) = \sum_{i=1}^{\infty} \mu(A_i)
\]
\end{defn}

\begin{defn}
A \emph{measure space} is a measurable space which has a measure defined on the $\sigma$-algebra of its measurable sets.
\end{defn}

\begin{defn}
The measure $\mu$ is called \emph{$\sigma$-finite} if $X$ is a countable union of measurable sets of finite measure.
\end{defn}

Let $\mu$ be a measure on a $\sigma$-algebra $\mathfrak{M}$, then
\begin{itemize}
	\item[(i)]	$\mu(\emptyset) = 0$.
	\item[(ii)]	$\mu(A_1 \cup \cdots \cup A_n) = \mu(A_1) + \ldots + \mu(A_n)$ if $A_1, \ldots, A_n$ are pairwise disjoint 											members of $\mathfrak{M}$.
	\item[(iii)]	$A \subseteq B$ implies $\mu(A) \leq \mu(B)$ for $A, B \in \mathfrak{M}$.
	\item[(iv)]	$\mu(A_n) \rightarrow \mu(A)$ as $n \rightarrow \infty$ if $A = \bigcup_{n=1}^{\infty}A_n,\;A_n \in 														\mathfrak{M}$, and $A_1 \subset A_2 \subset \cdots$
	\item[(v)]	$\mu(A_n) \rightarrow \mu(A)$ as $n \rightarrow \infty$ if $A = \bigcap_{n=1}^{\infty}A_n,\;A_n \in 														\mathfrak{M}$, and $A_1 \supset A_2 \supset \cdots$ and $\mu(A_1)$ is finite.
\end{itemize}

\begin{defn}
Let $(X, \tau)$ be a topological space; let $\mathfrak{B}$ denote the Borel $\sigma$-algebra on $X$, that is, the smallest $\sigma$-algebra on $X$ that contains all open sets $U \in \tau$.  Let $\mu$ be a measure on $\mathfrak{B}$.  Then $\mu$ is called a \emph{Borel measure}.
\end{defn}

\begin{defn}
The \emph{support} of a measure $\mu$ is defined to be the set of all points $x$ in $X$ for which every open neighbourhood of $x$ has positive measure:
\[
	\mbox{supp }(\mu) := \left \{ x \in X: x \in U \in \tau \Rightarrow \mu(U) > 0 \right \}
\]
\end{defn}

\begin{defn}
$L^1(X, \mu)$ is the set of equivalence classes of all real-valued functions on a topological space $X$ which have finite integral with respect to the measure $\mu$, where $f \sim g \mbox{ iff } \mu(\{x: f(x) \neq g(x)\}) = 0$.  It may be written as $L^1(X)$ when the measure is understood.  This space is equipped with the norm
\[
	||f||_1 = \int_X |f(x)| \, d\mu(x).
\]
\end{defn}


\begin{theo} [Egorov's theorem] \label{T:EgorovTheorem}
Given a sequence $(f_n)$ of real-valued functions on some measure space $(X, \mathfrak{M}, \mu)$ and a measurable set $A$ with $\mu(A) < \infty$ such that $(f_n)$ converges $\mu$-almost everywhere on $A$ to a limit function $f$, the following result holds: for every $\epsilon > 0$, there exists a measurable subset $B \subset A$ such that $\mu(B) < \epsilon$, and $(f_n)$ converges to $f$ uniformly on $A \setminus B$.
\end{theo}

In other words, pointwise convergence on $A$ implies uniform convergence everywhere except on some subset $B$ of fixed small measure.

\subsection{Graph Theory}

A \emph{graph} is a pair $G = (V, E)$ of sets such that $E \subseteq [V]^2$, where $[V]^2$ denotes the collection of all 2-element subsets of $V$, and $V \cap E = \emptyset$.  The elements of $V$ are called \emph{vertices} (or \emph{points}) of $G$ and the elements of $E$ are called \emph{edges}.  We will write $xy$ for an edge between two elements $x,y \in V$.  The number of edges at each vertex is called the \emph{degree} of the vertex.  The number of vertices of a graph is called its \emph{order}.  If we let $G = (V, E)$ and $G^\prime = (V^\prime, E^\prime)$ be two graphs, then $G$ and $G^\prime$ are \emph{isomorphic} if there exists a bijection $\Phi: V \rightarrow V^\prime$ with $xy \in E \Leftrightarrow \Phi(x)\Phi(y) \in E^\prime$ for all $x, y \in V$.  Such a map $\Phi$ is called and \emph{isomorphism} and we write $G \cong G^\prime$.  If $V^\prime \subseteq V$ and $E^\prime \subseteq E$, then $G^\prime$ is a \emph{subgraph} of $G$, written $G^\prime \subseteq G$.  If $G^\prime \subseteq G$ and $G^\prime$ contains all the edges $xy \in E$ with $x, y \in V^\prime$, then $G^\prime$ is an \emph{induced subgraph} of $G$.


\section{Product Topology vs. Box Topology}

Let $X$ and $Y$ be topological spaces, we shall call them \emph{coordinate spaces}.  Let $U \subseteq X$ and $V \subseteq Y$, with $U$ and $V$ open with respect to the topologies on $X$ and $Y$.  The family of all cartesian products $U \times V$ forms a base for a topology for $X \times Y$.  This topology is called the \emph{product topology} of $X \times Y$.

The functions $\pi_0: X \times Y \rightarrow X$ and $\pi_1: X \times Y \rightarrow Y$, which take $(x,y) \in X \times Y$ to $x \in X$ and $y \in Y$ respectively, are called \emph{projections} onto the coordinate spaces.  These functions are continuous because given any open set $U \subseteq X$, we have $\pi_0^{-1}(U) = U \times Y$ which is an open set in $X \times Y$.  The product topology is the coarsest topology (topology with the fewest open sets) for which all projections into their respective coordinate spaces are continuous.  

Now let us generalise this idea to the case of an arbitrary number of coordinate spaces.  Let $\{X_i: i \in I\}$ be a collection of topological spaces indexed by $I$ ($I$ may be finite, countable or uncountable), and let $X = \prod_i X_i$ be the cartesian product of these topological spaces.  A subbase for the product topology is formed by the collection of sets $\pi_i^{-1}(U)$ where $U \subseteq X_i$ is open.  A base for the product topology is the family of all finite intersections of these subbase elements.  

If our index set $I$ contains an  infinite number of elements, then as a base we have the family of all cartesian products of open sets from each coordinate space with all but finitely many factors equal to the entire space.  

A more na\"ive approach to the above generalisation leads to what is known as the box topology.  In this approach, our base consists of the family of sets formed by taking the product of infinitely many open subsets, one in each coordinate space.  Examples of topologies other than the box topology and the product topology on the cartesian product of a collection of topological vector spaces are given in \cite{CJK}.


It is obvious that for a finite number of coordinate spaces the box and product topologies agree.  It is also obvious that the box topology is finer than the product topology.  However, the most important reason why the product topology is chosen over the box topology is because many theorems about finite products hold for arbitrary products if the product topology is used, but not for the box topology \cite{JRM}.  That being said, the box topology is still useful for constructing counterexamples.  The box topology came before the product and was studied first.  It was not until Tychonoff that the product topology became the canonical topology for a cartesian product of an infinite number of topological spaces.

Let us see how these two topologies differ by some examples.  The following examples are taken from \cite{DS}.

\begin{egg}
In the product topology, the product of compact spaces is compact --- this is the famous Tychonoff product theorem.  This fails in the box topology.  Consider $\II^{\infty}$ --- the countable product of copies of the unit interval, $\II$.  If $A_0 = [0,1) \mbox{ and } A_1 = (0,1]$, then the collection of all open sets of the form $A_{\epsilon_1} \times A_{\epsilon_2} \times \cdots$, where $\epsilon_i = 0,1$, is an uncountable open cover of $\II^{\infty}$ with no proper subcover.  For if $A_{\epsilon_1} \times A_{\epsilon_2} \times \cdots$ is excluded from the cover, the point $(\epsilon_1,\epsilon_2,\ldots)$ is not covered.
\end{egg}

\begin{egg}
	In the product topology, the product of connected spaces is connected.  This is not true in the box topology.  Consider, for 		example, $\RR^{\infty}$ --- the countable product of real lines.  The set 
	\[
		A = \{(x_1,x_2,\ldots): {x_i} \mbox{ is a bounded sequence}\}
	\]
	is both open and closed in the box topology, and thus $\RR^{\infty}$ is not connected with respect to the box topology.
\end{egg}




The following theorem is taken from \cite{CJK}:

\begin{theo} \label{T:DiscontBoxProd}
Let $\II^\infty$ be equipped with the box topology.  Let $x, y \in \II^\infty$. Then $x$ and $y$ are in the same connected component of $\II^\infty$ if and only if the set $\{i: x_i \neq y_i\}$ has finite cardinality.  If the number of members of $\{i: x_i \neq y_i\}$ is infinite then there exists a $U \subseteq \II^\infty$ which is open and closed at the same time and for which $x \in U$ while $y \notin U$.
\end{theo}

\begin{rmk}
Here are some other interesting facts about the box topology:
\begin{itemize}
	\item Every basic open set in the product topology is in the box topology but not every box open set is in product topology.  Basic open sets are countable unions of intersections of finitely many sets of the form $\pi_i^{-1}(U)$ where $U \subseteq X_i$ is open.  On the other hand, box open sets can be arbitrary unions of intersections of infinitely many sets of the form $\pi_i^{-1}(U)$.
	\item	Parallelepipeds are an infinite-dimensional generalisation of the cuboid with sides parallel to the coordinate axes.  They are not open in the product topology but they are open in the box topology.
\end{itemize}
\end{rmk}


So, as follows from theorem \ref{T:DiscontBoxProd}, $\II^\infty$ with the box topology has uncountably many disconnected components, and because of this there are lots of continuous functions.  
This is a strange property of the box topology --- a property that we do not want, and this is one reason we will do away with the box topology.  The purpose for the brief appearance of the box topology here is that during the course of our research we put our hopes on it having a connection with the Lebesgue density theorem, but these hopes were dashed.

Next, we turn our attention to the issue at hand.  We give a brief introduction of the Lebesgue density theorem in finite dimension and the reason we have chosen to extend it.





\cleardoublepage

\chapter{Lebesgue Density Theorem}

The purpose of this work is to prove an extension of the theorem --- which goes by the name of this chapter --- to infinite-dimensional spaces and give examples of functions which satisfy it.  But first, seeing that there is no obvious connection between the Lebesgue density theorem and random graph theory, we show what role the density theorem plays in graph theory by going through a proof which uses it.  Then we state what the Lebesgue density theorem and examine it superficially.  Later on we attempt to explain why it is non-trivial to find an extension to infinite dimensions.

\section{Justification for Use}

We came across the Lebesgue density theorem trying to obtain an analogue of the proposition which follows (the proposition and some of the terms used here are taken from \cite{MP}, chapter 3).  The proposition may seem incomprehensible if we do not explain some of the notation first.  So let's get right to it.

We will use $\Gamma$ to denote a \emph{feasible} connected graph of order $k$ which means that the probability that some random geometric graph on $k$ vertices with radius $r$ is isomorphic to $\Gamma$ is strictly positive, for some $r > 0$.  
As an example, take the 2-dimensional case with the normal Euclidean distance.  The star graphs $S_1 \ldots S_6$ are feasible, but $S_7$ and above are not.  A star graph on $n$ vertices, $S_n$, is a graph with one vertex having degree $n-1$ and the remaining $n-1$ vertices each having degree $1$.  We will call a vertex of degree $1$ a \emph{leaf}.  If we assume that the radius from the centre vertex to each leaf is one, then for $S_7$, which has 6 leaves, there are two consecutive leaves for which the distance between them will be less than or equal to one and so we would need to have the edge between those two leaves added in order to make it a random geometric graph.  However, the additional edge leaves us with a graph which is not $S_7$, and thus the probability of getting a graph isomorphic to $S_7$ is zero.

Moving on, let $(x_1, x_2) \succeq (y_1, y_2)$ if and only if either $x_1 > y_1$ or $x_1 = y_1$ and $x_2 \geq y_2$, then $\succeq$ is called the lexicographic order on the plane $\RR^2$.  The generalisation to higher dimensions is obvious.  Now, given a finite set of points $\mathfrak{Y} \subset \RR^d$, let the first element of $\mathfrak{Y}$ according to the lexicographic ordering of $\RR^d$ be called the left-most point of $\mathfrak{Y}$ (LMP($\mathfrak{Y}$)).  Let $G_{n,A}$ be the number of unlabelled induced subgraphs of $G$, with radius $r_n$, isomorphic to $\Gamma$ for which the left-most point of the vertex set lies in $A \subseteq \RR^d$, and $E[G_{n,A}(\Gamma)]$ is the expectation of that number.  $\partial A$ denotes the boundary of $A$ and is defined as the intersection of the closure of $A$ with the closure of $A$s complement.  We take the Lebesgue measure of the boundary of $A$ to be zero.

Given a connected graph $\Gamma$ on $k$ vertices, and $A \subseteq \RR^d$, then for all finite $\mathfrak{Y} \subset \RR^d$:
\[
	h_\Gamma(\mathfrak{Y}) := 1_{\{G(\mathfrak{Y};1) \cong \Gamma\}}
\]
\[
	h_{\Gamma, n, A}(\mathfrak{Y}) := 1_{\{G(\mathfrak{Y};r_n) \cong \Gamma\} \cap \{LMP(\mathfrak{Y}) \in A\}}
\]
The first indicator function is 1 when the point set $\mathfrak{Y}$ with a radius of $r = 1$ forms a graph isomorphic to $\Gamma$.  For the second function, it is 1 when the point set $\mathfrak{Y}$ with a radius of $r = r_n$ forms an isomorphic graph with its left-most point coming from the set $A$.  It should be noted that the values of either above function is zero on a set $\mathfrak{Y}$ of less than $k$ vertices.

Finally, we set
\[
	\mu_{\Gamma, A} := k!^{-1} \int_A f(x)^k dx \int_{(\RR^d)^{k-1}} h_\Gamma(\{0,x_1, \dots, x_{k-1}\})d(x_1, \dots, x_{k-1})
\]
Points to note about this function $\mu_{\Gamma, A}$:
\begin{enumerate}
	\item $k!^{-1}$ makes sure that we count each graph only once since there are $k!$ ways to choose the same $k$ points.
	\item	$\int_A f(x)^k dx = 1$ when $A = \RR^d$.
	\item	$\int_{(\RR^d)^{k-1}} h_\Gamma(\{0,x_1, \dots, x_{k-1}\})d(x_1, \dots, x_{k-1})$ gives the proportion of the graphs 						isomorphic to $\Gamma$.
\end{enumerate}

Here, $f$ is some specified probability density function on $\RR^d$ which is bounded.  It should not be confused with the $f$ used elsewhere in this work to mean a probability density function on $\RR^\infty$.

\begin{prop}\label{T:SubgraphCount}
	Suppose that $\Gamma$ is a feasible connected graph of order $k \geq 2$, that $A \subseteq \RR^d$ is open with 									$\lambda^d(\partial A) = 0$, and that $\lim_{n \rightarrow \infty}(r_n) = 0$.  Then
	\[
		\lim_{n \rightarrow \infty} r_n^{-d(k-1)} n^{-k} E[G_{n,A}(\Gamma)] = \mu_{\Gamma,A} \qquad (*)
	\]
\end{prop}

\begin{proof}
	Clearly $E[G_{n,A}(\Gamma)] = {n \choose k} E[h_{\Gamma, n, A}(X_k)]$.  Hence,
	\begin{align}
		E[G_{n,A}(\Gamma)] = & {n \choose k} \int_{\RR^d} \ldots \int_{\RR^d} h_{\Gamma, n, A}(\{x_1, \ldots, x_k\})f(x_1)^k dx_k 					\ldots dx_1 \notag \\
		&+ {n \choose k} \int_{\RR^d} \ldots \int_{\RR^d} h_{\Gamma, n, A}(\{x_1, \ldots, x_k\}) \notag \\
		&\times \left ( \prod_{i=1}^k f(x_i) - f(x_1)^k \right ) \prod_{i=1}^k dx_i. \qquad (**) \notag
	\end{align}
	By the change of variables $x_i = x_1 + r_{n}y_{i}$ for $2 \leq i \leq k$, and $x_1 = x$, the first term on the right-hand 			side of (**) equals
	\[
		{n \choose k} r_n^{d(k-1)} \int_{\RR^d} \ldots \int_{\RR^d} h_{\Gamma, n, A}(\{x, x + r_{n}y_{2}, \ldots, x + r_{n}y_{k}\}) 		dy_k \ldots dy_2 f(x)^k dx.
	\]
	Since $A$ is open, for $x \in A$ the function $h_{\Gamma, n, A}(\{x, x + r_{n}y_{2}, \ldots, x + r_{n}y_{k}\})$ equals 					$h_\Gamma(\{0, y_2, \ldots, y_k\})$ for all large enough $n$, while for $x \notin A \cup \partial A$ it equals zero for all 		$n$.  Also, $h_{\Gamma, n, A}(\{x, x + r_{n}y_{2}, \ldots, x + r_{n}y_{k}\})$ is zero except for $(y_2, \ldots, y_k)$ in a 			bounded region of $(\RR^d)^{k-1}$, while $f(x)^k$ is integrable over $x \in \RR^d$ since $f$ is assumed bounded.  Therefore, 		by the dominated convergence theorem for integrals, the first term on the right-hand side of (**) is asymptotic to $n^k 					r_n^{d(k-1)} \mu_{\Gamma, A}$.
	
	On the other hand, the absolute value of the second term on the right-hand side of (**) multiplied by $n^{-k} r_n^{-d(k-1)}$ 		is bounded by $\int_{\RR^d} w_n(x_1) f(x_1) dx_1$, where we set
	\[
		w_n(x) := \int_{B(x;kr_n)} \ldots \int_{B(x;kr_n)} r_n^{-d(k-1)} \left | \prod_{i=2}^k f(x_i) - f(x)^{k-1} \right | dx_2 					\ldots dx_k.
	\]
	If $f$ is continuous at $x$, then clearly $w_n(x)$ tends to zero.  Even if $f$ is not almost everywhere continuous, we assert 	that $w_n(x)$ still tends to zero if $x$ is a Lebesgue point of $f$.  This is proved by induction of $k$; the inductive step 		is to bound the integrand by
	\[
		r_n^{-d(k-1)} \left ( |f(x_k) - f(x)| \prod_{i=2}^{k-1} f(x_i) \right ) + r_n^{-d(k-1)} \left | \prod_{i=2}^{k-1} f(x_i) - 			f(x)^{k-2} \right | f(x).
	\]
	The integral of the first expression over $B(x;kr_n)^{k-1}$ tends to zero by the definition of a Lebesgue point (and 						boundedness of $f$), while that of the second tends to zero by the inductive hypothesis.  Hence, by the Lebesgue density 				theorem and the dominated convergence theorem, $\int_{\RR^d} w_n(x_1) f(x_1) dx_1$ tends to zero, which proves the second 			equality in (*).
\end{proof}
So, as we see, the Lebesgue density theorem is of use in dispensing with any assumption of continuity on $f$.  For this reason we would like to have an analogue of it in the infinite-dimensional case.  Before we get to that, let us see what exactly this Lebesgue density theorem is and how it works in finite dimensions.  Then, once we have an understanding of the basics, we proceed to expand on that knowledge.

\section{The Classical Case for $\RR^d$}

To understand the Lebesgue density theorem, we first need to grasp the concept of a Lebesgue point.  Let us look at the definition.

\begin{defn}
	If $f \in L^1(\RR^d)$, any $x \in \RR^d$ for which it is true that
	\[
		\lim_{\epsilon \to 0} \frac{1}{\lambda^d(B_x(\epsilon))} \int_{B_x(\epsilon)} |f(y) - f(x)|d\lambda^d(y) = 0,
	\]
	where $\lambda^d$ is the $d$-dimensional Lebesgue measure, is called a \emph{Lebesgue point} of $f$.
\end{defn}

The definition of a Lebesgue point is simple and straightforward, but the Lebesgue density theorem is amazing and takes some time to get used to.

\begin{theo}[Lebesgue density theorem]
	If $f \in L^1(\RR^d)$, then almost every $x \in \RR^d$ is a Lebesgue point of $f$.
\end{theo}

The Lebesgue density theorem has a well-known proof which can be found in \cite{WR} p. 139.  Here we are going to work on it and transform it into some form that we can use.  We cannot use the proof as is for a couple reasons, the first is that it is proved for finite dimensions only, and another is that the tools used to prove it, such as the Hardy-Littlewood maximal function, may not survive the generalisation to infinite dimensions.  In fact, the Hardy-Littlewood maximal function has properties which are proved using the Vitali covering theorem which itself does not hold in the infinite-dimensional setting in which we are working.  This fact was proven by David Preiss \cite{DP} in 1979 and the result was strengthened by Jaroslav Ti\v ser \cite{JT1} in 2003.  However, the Lebesgue differentiation theorem, of which the Lebesgue density theorem is a special case, has been shown to hold for some class of Gaussian measures and all integrable functions provided that we change almost everywhere convergence to convergence in measure.  This is so even though the Vitali covering theorem fails in general for Gaussian measures \cite{JT2}.

\begin{rmk}
It is easy to see that the theorem holds true for locally integrable functions, but for our purpose in the case of $\RR^\infty$, local integrability is not defined.
\end{rmk}

Let us look at how the Lebesgue density theorem works for continuous functions.  All continuous functions are locally $L^1$ and so satisfy the theorem.  That is, every point has a neighbourhood such that the restriction of the function to this neighbourhood is $L^1$ (has a finite integral), and because the restriction to this neighbourhood is bounded and since our function is continuous, it is in $L^1$.  If $f$ is continuous, then the following proof (taken from \cite{NV}) shows that every member, $x$, of the domain of $f$ is a Lebesgue point.

\begin{theo} \label{T:LebAtContinuousPt}
	Given a measure, $\mu$, and a function, $f$ in $\RR^d$, which is continuous at one point, $x_0$.  Then $x_0$ is a Lebesgue point.
\end{theo}

\begin{proof}
	Since $f$ is continuous at $x_0$, then for all $\varepsilon > 0$, there exists a cube, $\Pi$, of strictly positive measure containing $x_0$ such that $\forall \, $y$ \, \in \, \Pi, |f(y) - f(x_0)| \leq \varepsilon$. Suppose $\mu(\Pi^{\prime}) > 0$ is such that $0 < \mu(\Pi^{\prime}) < \mu(\Pi)$, then:
	\[
		\frac{1}{\mu(\Pi^{\prime})} \int_{\Pi^{\prime}} |f(y) - f(x_0)| \, dy \leq \frac{1}{\mu(\Pi^{\prime})} \int_{\Pi^{\prime}} 				\varepsilon \, dy =	\varepsilon
	\]
	Therefore:
	\[
		\lim_{\mu(\Pi^{\prime}) \rightarrow 0} \frac{1}{\mu(\Pi^{\prime})} \int_{\Pi^{\prime}} |f(y) - f(x_0)| \, dy = 0
	\]
	and $x_0$ is a Lebesgue point of $f$.
\end{proof}
As a result of this simple proof, we are assured that every point of a continuous function is a Lebesgue point. Note that although $\Pi$ was chosen to be a cube, it could have been a ball or any Borel neighbourhood of $x$.



\begin{egg}
	The function \(
		f(x) = \left \{	
			\begin{array}{cl}
				0	&	\mbox{ if } x \in \RR \setminus \mathbb{Q} \\
				x	&	\mbox{ if } x \in \mathbb{Q}
			\end{array}
		\right.
	\)
	is continuous at only one point, namely $x = 0$, and in particular that point is a Lebesgue point of the function.
\end{egg}
	
At this moment, we should clarify a point which may have been missed.  A point of continuity of $f$ is a Lebesgue point of $f$, but a Lebesgue point of $f$ is not necessarily a point of continuity of $f$.  The following example demonstrates this fact.

\begin{egg}
	The function \(
		1_{\mathbb{Q}}(x) = \left \{
			\begin{array}{cl}
				0	&	\mbox{ if } x \notin \mathbb{Q} \\
				1	&	\mbox{ if } x \in \mathbb{Q}
			\end{array}
		\right.
	\)
	is nowhere continuous, however, every point $x \in \RR \setminus \mathbb{Q}$ is a Lebesgue point of $f$ even if every point is a point of discontinuity of $f$.
\end{egg}


As mentioned earlier, the Lebesgue density theorem is a special case of the Lebesgue differentiation theorem and we may sometimes prove results for the differentiation theorem which will still hold for the density theorem.  We will not stress the distinction between the two concepts, in fact, they are very similar.  The Lebesgue differentiation theorem states that:

\begin{theo}
For almost every $x \in \RR^d$:
\[
\lim_{\epsilon \to 0} \frac{1}{\lambda^d(B_x(\epsilon))} \int_{B_x(\epsilon)} f(y) \, d\lambda^d(y) = f(x)
\]
where $f: \RR^d \to \RR$ is Lebesgue-integrable, and $B_x(\epsilon)$ is the ball of radius $\epsilon$ around $x$.
\end{theo}

\cite{NV} demonstrates how easy it is to show that if $x \in \RR^d$ is a Lebesgue point of $f$, then
\[
\lim_{\epsilon \to 0} \frac{1}{\lambda^d(B_x(\epsilon))} \int_{B_x(\epsilon)} f(y) \,  d\lambda^d(y) = f(x).
\]
Let us assume $x \in \RR^d$ is a Lebesgue point of $f$.  Then for all $\epsilon > 0$, we have:
\begin{align}
\left | \frac{1}{\lambda^d(B_x(\epsilon))} \int_{B_x(\epsilon)} f(y) \, d\lambda^d(y) - f(x) \right | &= \left | \frac{1}{\lambda^d(B_x(\epsilon))} \int_{B_x(\epsilon)} (f(y) - f(x)) \, d\lambda^d(y) \right | \notag \\
	&\leq \frac{1}{\lambda^d(B_x(\epsilon))} \int_{B_x(\epsilon)} |f(y) - f(x)| \, d\lambda^d(y) \notag
\end{align}
Hence, from
\[
\lim_{\epsilon \to 0} \frac{1}{\lambda^d(B_x(\epsilon))} \int_{B_x(\epsilon)} |f(y) - f(x)| \, d\lambda^d(y) = 0
\]
we conclude that
\[
\lim_{\epsilon \to 0} \frac{1}{\lambda^d(B_x(\epsilon))} \int_{B_x(\epsilon)} f(y) \, d\lambda^d(y) = f(x).
\]

We have just scratched the surface of the Lebesgue density theorem.  Later on, we look at densities of functions in $\II^\infty$ and $\RR^\infty$ and we also look at a more general forms of the Lebesgue density theorem.  As stated earlier, the Lebesgue density theorem, as is, fails in infinite dimensions but we will attempt to remedy this by finding the correct version of the theorem which will work in infinite dimensions and in particular, $\RR^{\infty}$.  We will also detail the difficulties faced in extending it.

Before we do that, however, we need an interesting tool.  As seen in the proof of \ref{T:SubgraphCount}, it is of importance that we have a translation-invariant measure because we always rescale our points according to the left-most point, and if our measure is not translation-invariant then it is obvious that we would not be able to count the number of subgraphs using the same proof.  
This measure we are interested in should be on $\RR^\infty$ and we should be able to take projections down into any space of finite dimension $\RR^d$.  Once we find this measure, we will need a Lebesgue density theorem for it.  Next, we reveal the resulting candidate of our search for such a non-trivial translation-invariant measure on $\RR^\infty$.


\cleardoublepage

\chapter{``Lebesgue Measure'' on $\RR^\infty$}

Is it possible to measure the volume of an infinite-dimensional object?  If so, how do we do it?  And is it possible to do so in an easy and intuitive manner?  Can an infinite-dimensional object even have finite volume?  Of course it can, take as an example the infinite-dimensional cube of side 1, which we will denote by $C_1$.  Its volume is expected to be 1.  Now consider the parallelepiped (an infinite-dimensional generalisation of the cuboid with sides parallel to the coordinate axes) formed from $C_1$ by shortening one of its sides to $\frac{1}{2}$.  Its volume, also, is expected to be $\frac{1}{2}$.  Thus, from this example, we see that there exists uncountably many infinite-dimensional objects with finite volume.  The purpose of this chapter is to find a ``nice'', intuitive way to measure the volumes of such objects.  By nice we mean as simple as the Lebesgue measure for finite-dimensional objects.  Of course, we also want to go beyond well shaped objects such as cubes and parallelepipeds.  We need to make it as general as possible.

If a Lebesgue measure on $\RR^\infty$ already existed, then there would be no reason for this chapter.  And in the true sense of what a measure is, a non-trivial one does not exist and cannot exist.  Clearly, the trivial measure where we assign a `size' of zero to each set will work, but we require something with more substance.  Now, if we lessen our expectations slightly, we may be able to obtain what we desire.  We will settle for a measure which is not $\sigma$-finite.  Our setting requires a nice translation-invariant measure on $\RR^\infty$.  However, it is well known in functional analysis that the trivial measure is the only $\sigma$-finite translation-invariant Borel measure on an infinite-dimensional locally convex topological vector space \cite{HSY, YY}.  The proof of this statement can be found in \cite{YY} pp. 138 - 143.  As such, there is no analogue to the Lebesgue measure on a space of infinite dimensions.  Although this is well known, it does not stop mathematicians from coming up with measures on these spaces.  Richard Baker, amongst others \cite{TVY, YY}, has shown the existence of and constructed a non-trivial translation-invariant Borel measure which is almost as nice to work with as the Lebesgue measure, but on the infinite-dimensional space $\RR^\infty$.  It should be noted that this measure is not $\sigma$-finite as we mentioned earlier.  What does it mean when we say that this measure is almost as nice to work with as the Lebesgue measure?  As Baker puts it, it means that if $R = \prod_{i=1}^{\infty}(a_i,b_i)$ is any infinite-dimensional parallelepiped such that the ``volume'' $\prod_{i=1}^{\infty}(b_i-a_i)$ of $R$ is a non-negative real number, then 
\[
	\lambda^\infty(R) = \prod_{i=1}^{\infty}(b_i-a_i)
\]
where $\lambda^\infty$ is our so-called \emph{infinite-dimensional Lebesgue measure on $\RR^\infty$}.  This is as intuitive and as easy as it gets.  In this chapter we attempt to summarize (and reproduce some proofs) of Richard Baker's paper which goes by the same title as the chapter.  This is a survey of his paper and so no new results will be given in this chapter except for a lemma toward the end of it which tries to help us understand $\lambda^\infty$.  We will soon go through the construction of this measure in detail, but first, we provide a weaker proof than the one which states that there is only one $\sigma$-finite translation-invariant Borel measure on an infinite-dimensional space and it is the trivial measure.


\section{Impossibility Theorem}

One important thing to note is that the measure we are going to use, as simple and intuitive as it is, lacks a key property, that of being $\sigma$-finite.  This section attempts to show that it is actually impossible to have such a nice measure in $\RR^\infty$, not including the trivial measure.  To see that this measure is not $\sigma$-finite, let us refer back to our example of the cube.  However, instead of having a cube with sides of length 1, let it have sides of length 2 and denote it by $C_2$.  This cube exists in $\RR^\infty$.  The volume of this cube is obviously infinity, but the question is, can it be covered by a countable number of cubes of finite volume each?  The answer is no, for if we were to cover it with cubes of side length 1 (which have finite volume as noted earlier, that being 1), we would need $2^\infty$ of these cubes to cover $C_2$, and that is an uncountable number.  

\begin{theo}
	Let $(X, \| \cdot \|)$ be an infinite-dimensional, separable Banach space. Then the only locally finite (every point of the measure space has a neighbourhood of finite measure) and translation-invariant Borel measure $\mu$ on $X$ is the trivial measure, with $\mu(A) = 0$ for every measurable set $A \in X$.  
\end{theo}

\begin{proof}
Equip $X$ with a $\sigma$-finite, translation-invariant Borel measure $\mu$. 	Using $\sigma$-finiteness, choose a $\delta > 0$ such that the open ball $B_x(\delta)$ of radius $\delta$ around an arbitrary element $x \in X$ has finite	$\mu$-measure.  Since $X$ is infinite-dimensional, the ball $B_x(\delta)$ is non-compact, and so there exists a $\gamma > 0$ such that $B_x(\delta)$ cannot be covered with finitely many balls of radius $\gamma$.  By this fact, it is easy to construct an infinite sequence of points $x_1, x_2, \ldots, x_n, \ldots \in B_x(\delta)$ so that the open balls $B_{x_n}(\frac{\gamma}{2}), n \in N$, of radius $\frac{\gamma}{2}$, do not intersect, and are contained in $B_x(\delta)$.

By translation-invariance, all of the smaller balls $B_{x_n}(\frac{\gamma}{2})$ have the same measure and as the sum of these measures is finite, the smaller balls must all have $\mu$-measure zero.  Now, since $X$ is separable, it can be covered by a countable collection of balls of radius $\frac{\gamma}{2}$ and because each of these balls has $\mu$-measure zero, so must the whole space $X$, thus, $\mu$ is the trivial measure.
\end{proof}
Simply speaking, you can fit inside a ball, infinitely many smaller balls of equal size, which means by $\sigma$-additivity that each of these smaller balls has measure zero.  On the other hand, because the space is separable, you can cover it by countably many small balls of measure zero which means that it, too, has measure zero.  


\section{Construction of the Measure $\lambda^\infty$}

Here we get to the heart of the chapter and one of the main points of interest of this paper.  This construction of an analogue to the Lebesgue measure follows the same plan as that of the regular Lebesgue measure construction given in C. A. Rogers' ``Hausdorff Measures'' \cite{CAR}.  The idea behind this measure is as follows: cover our object (set) with as few infinite-dimensional parallelepipeds of finite volume as needed, then take the sum of the volumes of these parallelepipeds.  This gives a rough estimate of the volume of the object.  As there may be overlaps of these covering parallelepipeds, we take the smallest sum that we can get (which is obtained where the least overlap occurs).  Since it may not be possible to cover the object with totally disjoint parallelepipeds, the infimum of the above sum will suffice.

At this point, we would like to introduce some new notation and recall a few definitions, namely that of: positively separated (\ref{D:PositivelySeparated}), pre-measure (\ref{D:Premeasure}), outer measure (\ref{OuterMeasure}) and metric outer measure (\ref{D:MetricOuterMeasure}).  Then we state a few theorems without proof.  The proofs of these theorems and most of the terms defined here can be found in Chapter 1 of \cite{CAR}.

\begin{defn}
	Let $\mathfrak{R}$ be the family of all infinite-dimensional parallelepipeds $R~\in~\RR^{\infty}$ of the form
	\[
		R = \prod_{i=1}^{\infty}(a_i,b_i), -\infty < a_i \leq b_i < +\infty,
	\]
	such that $0 \leq \prod_{i=1}^{\infty}(b_i - a_i) < +\infty$ where the product converges in the normal sense.
\end{defn}

\begin{defn}
	Let $\xi$ be the function on $\mathfrak{R}$ defined by
	\[
		\xi(R) = \prod_{i=1}^{\infty}(b_i - a_i), \qquad R \in \mathfrak{R}.
	\]
\end{defn}

\begin{defn}
	Let $\lambda^\infty$ be the function defined on all subsets of $\RR^{\infty}$ by
	\[
		\lambda^\infty(E) = \inf_{\substack{R_j \in \mathfrak{R}\\\cup R_j \supseteq E}} \sum_{j=1}^{\infty} \xi(R_j), \qquad E \subseteq \RR^{\infty}.
	\]
\end{defn}

Let us agree that any infimum taken over an empty set of real numbers has the value of +$\infty$.  
	
\begin{theo} \label{T:TransInvBorelRinfty}
	The $\lambda^\infty$ as above is a translation invariant Borel measure on $\RR^{\infty}$ such that for all 
	$R = \prod_{i=1}^{\infty}(a_i,b_i) \in \mathfrak{R}$, we have
	\[
		\lambda^\infty(R) = \prod_{i=1}^{\infty}(b_i - a_i).
	\]
\end{theo}

Theorem \ref{T:TransInvBorelRinfty} is the main theorem of interest.  The purpose of the present chapter is to prove this theorem.

\begin{theo} [Method I]
	If $\xi$ is a pre-measure defined on a family $\mathfrak{C}$ of subsets of $X$, the function
	\[
		\mu(E) = \inf_{\substack{C_j \in \mathfrak{C}\\\cup C_j \supseteq E}} \sum_{j=1}^{\infty} \xi(C_j)
	\]
	is an outer measure on $X$.
\end{theo}

\begin{theo} \label{T:UniqueXsi}
	Let $X$ is an arbitrary non-empty set.  Let $\mu$ be an outer measure on $X$.  Let $\xi = \mu$, then $\xi$ is a pre-measure.  Let $\lambda$ be the outer measure constructed by Method I from the pre-measure $\xi$, then $\lambda$ coincides with $\mu$.
\end{theo}

\begin{theo} [Method II]
	If $\xi$ is a pre-measure defined on a family $\mathfrak{C}$ of subsets, in a metric space $X$ with metric $\rho$, the set 			function
	\[
		\mu(E) = \sup_{\delta > 0} \mu_\delta(E) = \lim_{\delta \rightarrow 0} \mu_\delta(E),
	\]
	where
	\[
		\mu_\delta(E) = \inf_{\substack{C_j \in \mathfrak{C}\\\mbox{diam}(C_j) \leq \delta\\\cup C_j \supseteq E}} \sum_{j=1}^{\infty} \xi(C_j)
	\]
	is an outer measure on $X$.  Here, $\mbox{diam}(C)$ is the diameter of $C$ with respect to $\rho$, $C \subseteq \RR^{\infty}$.
\end{theo}

\begin{theo}\label{T:MetricOuterMeasure}
	Let $\mu$ be an outer measure on $X$ (a metric space with metric $\rho$), constructed by Method II, from the pre-measure 				$\xi$.  Then $\mu$ is a metric outer measure on $X$.
\end{theo}

\begin{theo}\label{T:AllBorelMeasurable}
	If $\mu$ is a metric outer measure on a metric space $X$, then every Borel set in $X$ has a value uniquely defined by $\mu$.
\end{theo}

The plan for the remainder of the chapter is to construct the measure $\lambda^\infty$ on $\RR^\infty$ satisfying the properties that all Borel sets in $\RR^\infty$ have values uniquely defined by $\lambda^\infty$ and $\lambda^\infty$ is translation invariant on $\RR^\infty$.  We continue as follows:
\begin{itemize}
	\item By Method I and the definitions of $\mathfrak{R}$, $\xi$ and $\lambda^\infty$ we see that $\lambda^\infty$ is an outer measure on $\RR^\infty$.
	\item	Prove that $\forall R \in \mathfrak{R}, \xi(R) = \lambda^\infty(R)$.
	\item	Prove that $\forall R \in \mathfrak{R}, \nu(R) = \xi(R)$, where $\nu$ is an outer measure constructed by Method II.
	\item	Show that $\lambda^\infty(E) = \nu(E)$ for all $E \subseteq \RR^\infty$, and thus $\lambda^\infty$ is a metric outer measure by \ref{T:MetricOuterMeasure} above.
	\item	By \ref{T:AllBorelMeasurable} above, every Borel set in $\RR^\infty$ has a value unique defined by $\lambda^\infty$.
	\item	Finally, show that $\lambda^\infty$ is translation-invariant on $\RR^\infty$.
\end{itemize}

\begin{theo} \label{T:TauLeqLambda}
	Let $I = \prod_{i=1}^\infty[a_i, b_i],\; - \infty < a_i \leq b_i < + \infty$, be an infinite-dimensional compact 								parallelepiped in $\RR^\infty$ such that $0 \leq \prod_{i=1}^\infty(b_i - a_i) < + \infty$, then
	\[
		\prod_{i=1}^\infty(b_i - a_i) \leq \lambda^\infty(I).
	\]
\end{theo}





\begin{proof}
	If $I \subseteq \bigcup_{j=1}^\infty R_j$, where $R_j \in \mathfrak{R}$, then it is enough to show that
	\[
		\xi(I) \leq \sum_{j=1}^\infty\xi(R_j) \qquad (*)
	\]
	where $\xi(I) = \prod_{i=1}^\infty(b_i - a_i)$.  Assume the strict inequalities $0 < \xi(I)$ and $\sum_{j=1}^\infty\xi(R_j)~<~\infty$ hold as the proof is trivial in those cases where they are equal.  For all $d, j \geq 1$, let us introduce the following notation
	\[
		R_j = \prod_{i=1}^\infty(a_{ij}, b_{ij}), \qquad R_{dj} = \prod_{i=1}^d (a_{ij}, b_{ij}) \times \prod_{i=d+1}^\infty \RR.
	\]
	
	\noindent Choose $\epsilon > 0$.
	
	\noindent If we assume $\xi(I)$ to be finite, then there must exist a $d \in \NN$ such that:
	\begin{itemize}
		\item[(1)]	$\prod_{i=d+1}^\infty (b_i - a_i) < 1 + \epsilon$ (see lemma \ref{L:SidesApproachOne}).
	\end{itemize}
	
	\noindent Let $\mathfrak{F}$ be the family of sets containing parallelepipeds $R_{j}$ satisfying (1) above and either:
	\begin{itemize}
		\item[(2)]	$\prod_{i=d+1}^\infty(b_{ij} - a_{ij}) > 1 - \epsilon$; or
		\item[(3)]	$\prod_{i=1}^{d}(b_{ij} - a_{ij}) < \frac{\epsilon}{2^j}$.
	\end{itemize}
	  
	\noindent Since $0 < \xi(I) < + \infty$ and $0 \leq \xi(R_j) < + \infty$, $\mathfrak{F}$ clearly covers $I$ and as 						parallelepipeds in $\mathfrak{F}$ are open and $I$ is compact, there exists a finite subfamily $\{R_{d_{p}j_{p}} | 1 \leq p \leq 		k\}$ of $\mathfrak{F}$ that covers $I$.
	
	\noindent Choose $d > max\{d_1, \ldots, d_k\}$ and for $1 \leq p \leq k$, define
	\[
		I_d = \prod_{i=1}^d [a_i, b_i], \qquad S_{dp} = \prod_{i=1}^{d_p}(a_{ij_p}, b_{ij_p}) \times \prod_{i=d_p+1}^d [a_i, b_i].
	\]
	
	\noindent It is easy to see that $I \subseteq \bigcup_{p=1}^k R_{d_{p}j_{p}}$ and thus $I_d \subseteq \bigcup_{p=1}^k S_{dp}$.
	
	\noindent Let $\lambda^d$ be the usual Lebesgue measure on $\RR^d$, then
	\begin{align}
		\prod_{i=1}^d (b_i - a_i) = \lambda^d (I_d) &\leq \sum_{p=1}^k \lambda^d (S_{dp}) \notag \\
		&= \sum_{p=1}^k \left \{\prod_{i=1}^{d_p}(b_{ij_p} - a_{ij_p}) \cdot \prod_{i=d_p+1}^d (b_i - a_i) \right \}. \notag
	\end{align}
	
	\noindent Taking the limit as $d \rightarrow \infty$, we get
	\begin{align}
		\xi(I) &\leq \sum_{p=1}^k \left \{\prod_{i=1}^{d_p}(b_{ij_p} - a_{ij_p}) \cdot \prod_{i=d_p+1}^\infty (b_i - a_i) \right \} 						\notag \\
		&\leq (1 + \epsilon)\sum_{p=1}^k \prod_{i=1}^{d_p}(b_{ij_p} - a_{ij_p}) \mbox{ (by (1)) } \notag \\
		&\leq (1 + \epsilon)\sum_{p=1}^k {\,}^{{\!}^{\prime}} \prod_{i=1}^{d_p}(b_{ij_p} - a_{ij_p}) + 
				(1 + \epsilon)\sum_{p=1}^k {\,}^{{\!}^{\prime\prime}} \prod_{i=1}^{d_p}(b_{ij_p} - a_{ij_p}) \mbox{ (by (2) and (3)) } 					\notag
	\end{align}
	where $\Sigma^{{\,}^\prime}$ is the sum over those $p$ for which $\prod_{i=d_p+1}^\infty (b_{ij_p} - a_{ij_p}) > 1 - 						\epsilon$ 
	and $\Sigma^{{\,}^{\prime \prime}}$ is the sum over those $p$ for which 
	$\prod_{i=1}^{d_p}(b_{ij_p} - a_{ij_p}) < \frac{\epsilon}{2^{j_p}}$.  It follows that
	\begin{align}
		\xi(I) &\leq \left ( \frac{1 + \epsilon}{1 - \epsilon} \right ) \sum_{p=1}^k {\,}^{{\!}^{\prime}} \xi(R_{j_p}) + 
				(1 + \epsilon)\epsilon \sum_{p=1}^k {\,}^{{\!}^{\prime\prime}} \frac{1}{2^{j_p}} \mbox{ (by (3), (4)) } \notag \\
		&\leq \left ( \frac{1 + \epsilon}{1 - \epsilon} \right ) \sum_{j=1}^\infty \xi(R_j) + (1 + \epsilon)\epsilon \notag.
	\end{align}
	As $\epsilon \to 0$, (*) holds.
\end{proof}
\begin{theo} \label{T:TauEqLambda}
	For every $R \in \mathfrak{R}$, $\xi(R) = \lambda^\infty(R)$.
\end{theo}

\begin{proof}
Let $R \in \mathfrak{R}$.  Clearly $\lambda^\infty(R) \leq \xi(R)$, and we may assume that $\xi(R)~>~0$.  Let $\epsilon > 0$.  There exists a compact parallelepiped $I = \prod_{i=1}^\infty[a_i, b_i] \subseteq R$ such that $\prod_{i=1}^\infty(b_i~-~a_i)~=~(1~-~\epsilon)\xi(R)$.  By \ref{T:TauLeqLambda}, $\prod_{i=1}^\infty(b_i - a_i) \leq \lambda^\infty(I)$, thus \\$(1~-~\epsilon)\xi(R)~\leq~\lambda^\infty(R)$.  As $\epsilon \to 0$, we obtain the desired result.
\end{proof}
\begin{theo} \label{T:OuterMeasPreMeas}
	Let $\nu$ be the outer measure on $\RR^{\infty}$ constructed from the pair $\xi, \rho$ by Method II.  Then for all $R \in 			\mathfrak{R}$, $\nu(R) = \xi(R)$.
\end{theo}

\begin{proof}
	Let $R \in \mathfrak{R}$.  By \ref{T:TauEqLambda}, $\lambda^\infty(R) = \xi(R)$.  Clearly $\lambda^\infty(R) = \nu(R)$, and it is sufficient to show that
	\[
		\nu(R) \leq \xi(R). \qquad (**)
	\]
	Let $R = \prod_{i=1}^{\infty}(a_i,b_i), -\infty < a_i \leq b_i < +\infty$.  Assume $R \neq \emptyset$ as the proof is 				trivial if it is.  Thus, for all $i, a_i < b_i$.  Now, if $\xi(R) = 0$, then for all $d$, we have 											$\prod_{i=d+1}^{\infty} (b_i - a_i) = 0$.  However, if $\xi(R) > 0$, then $\lim_{d \rightarrow \infty} 												\prod_{i=d+1}^{\infty} (b_i - a_i) = 1$ (see lemma \ref{L:SidesApproachOne}).  Choose $\delta, \epsilon > 0$.  Then there exists an $d$ such that $\prod_{i=d+1}^{\infty} (b_i - a_i) < 1 + \epsilon$ and $\sum_{i=d+1}^{\infty} 2^{-i} < \frac{\delta}{2}$.  Define $R_d = \prod_{i=1}^{d} (a_i,b_i)$.  For $x = (x_i)_{i=1}^d, y = (y_i)_{i=1}^d \in \RR^d$, let $\rho_d(x,y)$ be defined as
	\[
		\rho_d(x,y) = \sum_{i=1}^{d} 2^{-i} \frac{|x_i - y_i|}{1 + |x_i - y_i|}.
	\]
	Cover the parallelepiped $R_d$ by parallelepipeds $R_{d1}, \ldots, R_{dm}$ in $\RR^d$ such that
	\begin{itemize}
		\item[(a)]	For $1 \leq j \leq m$, $R_{dj} = \prod_{i=1}^{d} (a_{ij},b_{ij}), , -\infty < a_{ij} \leq b_{ij} < +\infty$.
		\item[(b)]	For $1 \leq j \leq m$, $\sup_{x, y \in R_{dj}} ||x-y|| < \frac{\delta}{2d}$, where $||\cdot||$ is the Euclidean 			norm on $\RR^d$.
		\item[(c)]	$\sum_{j=1}^{m} \lambda^d (R_{dj}) < \prod_{i=1}^{d} (b_i - a_i) + \epsilon$.
	\end{itemize}
	For $1 \leq j \leq m$, define $R_j = R_{dj} \times \prod_{i=d+1}^{\infty} (a_i,b_i)$.  Let $x = (x_i), y = (y_i) \in R_j$, and define $x^{(d)} = (x_i)_{i=1}^d, y^{(d)} = (y_i)_{i=1}^d$, then, by (b), we have
	\begin{align}
		\rho_d(x,y) &= \rho_d(x^{(d)},y^{(d)}) + \sum_{i=d+1}^{\infty} 2^{-i} \frac{|x_i - y_i|}{1 + |x_i - y_i|} \notag \\
		&< \rho_d(x^{(d)},y^{(d)}) + \frac{\delta}{2} \notag \\
		&\leq d||x^{(d)}-y^{(d)}|| + \frac{\delta}{2} \notag \\
		&< \delta. \notag
	\end{align}
	Hence, for all $1 \leq j \leq m$, we have $\mbox{diam}(R_j) < \delta$.  It is clear that $R \subseteq \cup_{j=1}^m R_j$, hence by definition, $\nu_{\delta}(R) \leq \sum_{j=1}^{m} \xi(R_j)$.  By (c),
	\begin{align}
		\sum_{j=1}^{m} \xi(R_j) &= \sum_{j=1}^{m} \left \{ \lambda^d(R_{dj}) \cdot \prod_{i=d+1}^{\infty} (b_i - a_i) \right \} \notag \\
		&\leq \left \{ \prod_{i=1}^{d} (b_i - a_i) + \epsilon \right \} \cdot \prod_{i=d+1}^{\infty} (b_i - a_i) \notag \\
		&= \xi(R) + \epsilon \prod_{i=d+1}^{\infty} (b_i - a_i) \notag \\
		&< \xi(R) + \epsilon(1 + \epsilon). \notag
	\end{align}
	As $\epsilon \to 0$, we get $\nu_{\delta}(R) \leq \xi(R)$.  But $\delta \to 0$ also, hence $\nu(R) \leq \xi(R)$.  Therefore, (**) holds.
\end{proof}
\begin{theo}
	For all $E \subseteq \RR^{\infty}$, $\lambda^\infty(E) = \nu(E)$.  Hence, by Theorem \ref{T:MetricOuterMeasure}, $\lambda^\infty$ is a metric outer measure on $\RR^{\infty}$.
\end{theo}

\begin{proof}
	For $E \subseteq \RR^{\infty}$, we have $\lambda^\infty(E) \leq \nu(E)$, hence it suffices to prove that
	\[
		\nu(E) \leq \lambda^\infty(E), \qquad E \subseteq \RR^{\infty}. \qquad (***)
	\]
	Fix $E \subseteq \RR^{\infty}$. By Theorem \ref{T:UniqueXsi}, we have
	\[
		\nu(E) = \inf_{\substack{C_j \subseteq \RR^{\infty}\\\cup C_j \supseteq E}} \sum_{j=1}^{\infty} \nu(C_j).
	\]
	Hence we see that
	\[
		\nu(E) \leq \inf_{\substack{R_j \in \mathfrak{R}\\\cup R_j \supseteq E}} \sum_{j=1}^{\infty} \nu(R_j).
	\]
	For all $R \in \mathfrak{R}$, Theorem \ref{T:OuterMeasPreMeas} implies that $\nu(R) = \xi(R)$, therefore we have $\nu(E)~\leq~\lambda^\infty(E)$.  This proves (***).
\end{proof}
\begin{theo}
	The outer measure $\lambda^\infty$ is translation invariant on $\RR^\infty$.
\end{theo}

\begin{proof}
	This follows from the fact that if $R = \prod_{i=1}^{\infty}(a_i,b_i) \in \mathfrak{R}$ and $x \in \RR^\infty$, then $R + x 		\in \mathfrak{R}$ and $\xi(R + x) = \xi(R)$.
\end{proof}
To illustrate the properties of the measure $\lambda^\infty$ we prove the following lemma which we have already referred to and will be using again later on:

\begin{lem} \label{L:SidesApproachOne}
Let $\lambda^\infty(\Pi) = c$ with $0 < c < \infty$ where $\Pi = \prod_{i=1}^\infty [a_i, b_i]$ and each interval has positive length $l_i = b_i - a_i$.  Then $(l_i)_{i \to \infty} \to 1$.
\end{lem}

\begin{proof} 
If $\lambda^\infty(\Pi) = c$ then
\[
	\lambda^\infty(\Pi) = \lim_{d \to \infty} \prod_{i=1}^d l_i = c. \qquad (****)
\]
Taking $\log$s of both sides of (****) we get 
\begin{align}
	\lim_{d \to \infty} &\sum_{i=1}^d \log l_i \to \log c \notag \\
	&\Rightarrow \sum_{i=1}^\infty \log l_i < \infty \notag \\
	&\Rightarrow \log l_i \to 0 \notag \\
	&\Rightarrow l_i \to 1 \notag
\end{align}
Thus, $\forall \epsilon > 0, \exists N \mbox{ such that } \forall n \geq N, 1 - \epsilon < l_i < 1 + \epsilon$.
\end{proof}
It is interesting to note that lemma \ref{L:SidesApproachOne} implies that if we go to sufficiently high dimensions, a parallelepiped of finite volume looks more and more like the unit cube.  There are cases where this product may not be defined, like in
\[
\Pi = [0, \frac{1}{2}] \times [0,2] \times [0, \frac{1}{2}] \times [0,2] \times \ldots,
\]
but this does not concern us for we are only interested in ``nice'' parallelepipeds which are well behaved as in they have finite measure.

In trying to understand $\lambda^\infty$ let us look at the following example:

\begin{egg}
Denote the ``boundary'' of $\II^\infty$ by the following set:
\[
\partial \II^\infty = \{x \in \II^\infty: \exists i, x_i \in \{0,1\}\}.
\]
and let
\[
\partial_j \II^\infty = \{x \in \II^\infty: x_j = 0 \mbox{ or } 1\}.
\]
Thus, we can write $\partial \II^\infty = \cup_{j=1}^\infty \partial_j \II^\infty$.  Let $\epsilon > 0$ and write
\[
\epsilon = \sum_{j=1}^\infty 2^{-j} \cdot \epsilon.
\]
Let
\begin{align}
C_j = &[0,1] \times \ldots \times [0,1] \times [-2^{-j-2}\epsilon, 2^{-j-2}\epsilon] \times [0,1] \times \ldots \notag \\
&\cup [0,1] \times \ldots \times [0,1] \times [1-2^{-j-2}\epsilon, 1+2^{-j-2}\epsilon] \times [0,1] \times \ldots \notag
\end{align}
Since $\partial \II^\infty$ can be covered with countably many sets of volume zero, the ``boundary'' of $\II^\infty$ has measure zero.
\end{egg}

Having gone through the construction of this measure which was a bit technical, it prepared us for the next chapter which presents a paper by Jean Dieudonn\'e.  The chapter highlights the non-triviality of extending the Lebesgue density theorem to the infinite-dimensional case.
\cleardoublepage


\chapter{Dieudonn\'e's Example}\label{C:Dieudonne}

In this chapter we survey a paper by Jean Dieudonn\'e, one of the founding members of Bourbaki and a major contributor to the field of functional analysis.  The paper \cite{JD} is not famous and is only available in French.  Here, we reproduce it in English and attempt to explain its importance, both in our situation and more generally.  In brief, it explains the construction of a set which shows that there is no straightforward generalisation to infinite dimensions of the Lebesgue density theorem for finite dimensions.  If we put forward the most natural version of the extension, it is wrong and it will not work.


Before we get down to the details, let us set the foundation on which we will build.  First, the space we will work in is the Hilbert Cube, $\II^{\infty}$.  It is easier to work with because it has some desirable properties which are noticeably absent from $\RR^{\infty}$.  Later on we will undertake the task of extending our results from $\II^{\infty}$ to $\RR^{\infty}$.  Now, let $x \in \II^{\infty}$ be written as $(x^{\prime}, x^{\prime\prime})$ where $x^{\prime} \in \II^d$ and $x^{\prime\prime} \in \II^{\NN \setminus [1,2,\ldots,d]}$ and denote by $f_d$ the function obtained by integrating along the tail.  That is,
\[
	f_d(x^{\prime}) := \int_{\II^{\NN \setminus [1,2,\ldots,d]}} f(x^{\prime}, x^{\prime\prime}) d\lambda^{\NN \setminus [1,2,\ldots,d]}(x^{\prime\prime}).
\]
Thus, Jessen's theorem states that $f_d \rightarrow f$ almost everywhere as $d \rightarrow \infty$.

So a natural question to ask to make Jessen's theorem more general is: does it still work if we choose arbitrary finite subsets of the coordinate spaces?  To be more precise, given $J \in F$ where $F$ is the set of all finite subsets of $\NN$ ordered by inclusion, does  $f_J \rightarrow f$ almost everywhere as the cardinality of $J$ increases? 
The answer, unintuitively, is no, and we shall see why.

Also, as a consequence of Dieudonn\'e's example and to our detriment, if we have a bounded, measurable (thus integrable) function $f: \II^\infty \rightarrow \RR$, and the functions $f_d$ obtained by integrating along the tail.  Then, by looking at \[
	\lim_{d \rightarrow \infty} \left ( \lim_{n \rightarrow \infty} \frac{1}{\lambda^d(\Pi_n)} \int_{\Pi_n} 													|f_d(x^{\prime})| \, d \lambda^d \right )
\]
where $\lambda^d(\Pi_n) \rightarrow 0 \mbox{ as } n \rightarrow \infty$, the inner limit converges for fixed $d$ by the Lebesgue differentiation theorem, and the outer limit converges by Jessen's theorem.  However, because of Dieudonn\'e's example, we cannot send the outer and inner limits to infinity at the same time because we will not obtain $f(x)$.  It must be done in the correct order, fix a $d$, then send $n$ to infinity, then send $d$ to infinity.  Unfortunately, we are on the bad side of Dieudonn\'e's example.  Thus, the most general form of the Lebesgue density theorem does not hold in $\II^\infty$ and thus there is no hope for an extension to $\RR^\infty$.

\section{Dieudonn\'e's Example}

For the next two sections, to avoid excessive clutter, let us agree to denote all finite-dimensional or infinite-dimensional Lebesgue measures with $\mu$.  

Let us say that $F$ is a countable family of finite subsets of $\NN$ ordered by inclusion, thus $F$ is an ideal (see definition \ref{D:Ideal}).  It should be noted that the union of two finite sets is finite.  The theorem of Jessen would lead us to think that for almost every $x \in \II^\infty$, $f_J(x) \rightarrow f(x)$ according to the ideal $F$ (that is to say, for each $x \in \II^\infty$ not belonging to a set of measure zero and for each $\epsilon > 0$ there corresponds a finite subset $J_0(x, \epsilon)$ such that for all finite sets $J \supset J_0$, $|f_J(x) - f(x)| \leq \epsilon)$.
The goal of Dieudonn\'e's paper is to prove that the conjecture is inaccurate by building an example where this property fails.

Since the ideal $F$ has a countable basis, the reasoning which shows the Egorov theorem (see theorem \ref{T:EgorovTheorem}) for the almost everywhere convergence of a sequence of functions still applies when its almost everywhere convergence according to $F$ of a family $(g_J)$ of measurable functions for all $J \in F$.  More precisely, if such a family $(g_J) \rightarrow g$ almost everywhere according to the sets of the ideal $F$, for each $\delta > 0$, there exists a set $H \subset \II^\infty$, with $\mu(H) > 1 - \delta$, such that in $H$ the family $(g_J) \rightarrow g$ uniformly.  

In the example which Dieudonn\'e constructed, the function $f = 1_A$ is the characteristic (indicator) function of a set $A$ and $\mu(A) < \frac{1}{8}$.  For all $J \in F$, we will denote by $\phi_J(x) = \sup_{K \in J} f_K(x)$.  
If $f_J(x) \rightarrow f(x)$ almost everywhere it follows that there would exist a subset $J_0 \in F$ such that, for $J \supset J_0$, $\mu(\{x \in \II^\infty: \phi_J(x) > \frac{1}{2}\}) < \frac{1}{4}$.  However, the set $A$ will be such that there would exist an increasing sequence $J_n$ of finite subsets of $\NN$, such that for all $n$, $\mu(\{x \in \II^\infty: \phi_{J_n}(x) > \frac{1}{2}\}) > \frac{7}{16}$.  That will prove that $f_J(x) \nrightarrow f(x)$ almost everywhere according to F.

\subsection{Notation}

\noindent The sets $J_n = [1, q_n]$.\\

\noindent We will divide $J_n$ and $J_{n+1}$ into $h_n$ intervals labelled $J_{n,1}, J_{n,2}, \ldots, J_{n,h_n}$ and let $p_{n,r}$ denote the number of elements of $J_{n,r}$.  Thus, $q_{n+1} - q_n = p_{n,1} + p_{n,2} + \dots + p_{n,h_n}$.  The numbers $p_{n,r}$ and $h_n$ will be determined by induction and at the same time we define $A$ to be the union of pairwise disjoint sets $A_{n,r}$.\\

\noindent Denote by $k_n$ the total number of intervals $J_{m,r}$ $(m < n)$ which we divide $J_n$ into.  So $k_{n+1} = k_n + h_n$.\\

\noindent We have a decreasing sequence of positive numbers $(a_n)$ satisfying the following conditions:
\begin{itemize}
	\item[(a)]	$\sum_{n=1}^{\infty} a_n < \frac{1}{8}$ converges.
	\item[(b)]	$(a_n \log \frac{1}{a_n})$ diverges and each term $(a_n \log \frac{1}{a_n}) < \frac{1}{4}$.
\end{itemize}

\begin{egg}
	Taking
	\[
		a_n = \frac{c}{n(\log n)^2}
	\]
	for small enough $c$ will satisfy these conditions.
\end{egg}

\noindent We will denote $a_{k_n + r}$ by $a_{n,r}$ to simplify the notation.\\

\noindent Suppose that the $k_n$ sets $A_{m,r}$ $(m < n)$ form $\bar{A}_{m,r} \times I^{J_n^\prime}$ where $\bar{A}_{m,r} \subseteq I^{J_n}$.  Moreover, suppose $\mu(A_{m,r}) < a_{m,r}$ for the $k_n$ sets.  Let $\bar{B}_n = I^{J_n} \setminus \bigcup_1^{k_n} \bar{A}_{m,r}$, then one has $\mu_{J_n}(\bar{B}_n) > \frac{7}{8}$.

\subsection{Basis}

\noindent Let us start by defining the number $p_{n,1}$ and the set $A_{n,1}$.\\

\noindent So let $K_{n,1} = J_n \cup J_{n,1}$.\\

\noindent Take $A_{n,1} = \bar{B}_n \times \bar{C}_{n,1} \times I^{K_{n,1}^\prime}$, where $\bar{C}_{n,1} \subseteq I^{J_{n,1}}$ which we will define.\\

\noindent We will take $\bar{C}_{n,1} = \prod_1^{p_{n,1}} T_j$ (each $T_j$ taken in a interval $I_j$ of $I^{J_{n,1}}$ of length
\[
	1 - \frac{1}{p_{n,1}} \log \frac{1}{a_{n,1}}).
\]
Thus
\[
	\mu(\bar{C}_{n,1}) = \left (1 - \frac{1}{p_{n,1}} \log \frac{1}{a_{n,1}} \right ) ^ {p_{n,1}},
\]
$\mu(\bar{C}_{n,1}) \rightarrow 1$ as $p_{n,1} \rightarrow \infty$.\\

\noindent Take $p_{n,1}$ large enough so that $\frac{1}{2} a_{n,1} \leq \mu(\bar{C}_{n,1}) \leq a_{n,1}$ and
\[
	1 - \frac{1}{p_{n,1}} \log \frac{1}{a_{n,1}} > \frac{1}{2}.
\]

\noindent For each index $j \in J_{n,1}$, let $S_j= \prod_{\substack{i \in J_{n,1}\\i \neq j}} T_i$; $A_{n,1} = (\bar{B}_n \times S_j) \times (T_j \times I^{K_{n,1}^\prime})$.\\

\noindent For $J = J_n \cup (J_{n,1} - \{j\})$, we have $f_J(x) > \frac{1}{2}$ for all $x \in (\bar{B}_n \times S_j) \times (I_j \times I^{K_{n,1}^\prime})$.\\

\noindent So $D_{n,1}$ is the union of the $p_{n,1}$ sets for each $j \in J_{n,1}$.  It is clear that $\phi_{J_n}(x) > \frac{1}{2}$ in the set $D_{n,1}$.\\

\noindent We can write $D_{n,1} = \bar{B}_n \times \bar{D}_{n,1} \times I^{K_{n,1}^\prime}$, where $\bar{D}_{n,1} = \cup_1^{p_{n,1}} (S_j \times I_j)$.  It follows immediately that the measure of $\bar{D}_{n,1}$ (in $I^{J_{n,1}}$) is
\[
	\mu(\bar{D}_{n,1}) = \delta_{n,1} = \left (1 - \frac{1}{p_{n,1}} \log \frac{1}{a_{n,1}} \right ) ^ {p_{n,1}} + \left (\log \frac{1}{a_{n,1}} \right ) \left (1 - \frac{1}{p_{n,1}} \log \frac{1}{a_{n,1}} \right ) ^ {p_{n,1} - 1}.
\]

\noindent It is clear that $p_{n,1}$ can be large enough so that $\frac{1}{2} a_{n,1} \log \frac{1}{a_{n,1}} \leq \delta_{n,1} \leq 2 a_{n,1} \log \frac{1}{a_{n,1}}$.\\

\noindent We let $\bar{E}_{n,1} = I^{J_{n,1}} \setminus \bar{D}_{n,1}$ which has measure $1 - \delta_{n,1} > \frac{1}{2}$.

\subsection{Inductive Step}

\noindent Suppose now the sets $J_{n,1}, \ldots, J_{n,r}$ have all been defined and in each product $I^{J_{n,s}}$ for $s~\leq~r$, the two sets $\bar{C}_{n,s}$ and $\bar{D}_{n,s}$ are such that $\bar{C}_{n,s} \subset \bar{D}_{n,s}$ and that we have the following:
\begin{itemize}
	\item [(1)]	$H_{n,s}$ is the union of $J_{n,1}, \ldots, J_{n,s}$; $K_{n,s}$ is the union of $J_n$ and $H_{n,s}$.  Define by 								recurrence the sets $\bar{F}_{n,s}$ and $\bar{E}_{n,s}$ in $I^{H_{n,s}}$ as being the complement in the product, 								and as $\bar{F}_{n,1} = \bar{D}_{n,1}$ and
							\[
								\bar{F}_{n,s} = \left (\bar{F}_{n,s-1} \times I^{J_{n,s}} \right ) \cup \left (\bar{E}_{n,s-1} \times 													\bar{D}_{n,s} \right )
							\]
							then $A_{n,s} = \bar{B}_{n} \times \bar{E}_{n,s-1} \times \bar{C}_{n,s} \times I^{K_{n,s}^\prime}$ with 									$\mu(A_{n,s}) \leq a_{n,s}$.
	\item [(2)]	$\frac{1}{2} a_{n,s} \log \frac{1}{a_{n,s}} \leq \mu(\bar{E}_{n,s-1} \times \bar{D}_{n,s}) \leq 2 a_{n,s} \log \frac{1}{a_{n,s}}$.\\Moreover, $(\delta_{n,1} + \delta_{n,2} + \dots + \delta_{n,r}) \leq \frac{1}{2}$ in $\bar{F}_{n,r}$.
	\item [(3)]	For $D_{n,s} = \bar{B}_{n} \times \bar{E}_{n,s-1} \times \bar{D}_{n,s} \times I^{K_{n,s}^\prime}$, we have 											$\phi_{J_n}(x) > \frac{1}{2}$ on the set $D_{n,s}$.
\end{itemize}

\noindent So, $\mu(\bar{E}_{n,r}) = \beta_r = 1 - (\delta_{n,1} + \delta_{n,2} + \dots + \delta_{n,r}) \geq \frac{1}{2}$.\\

\noindent We take $\bar{C}_{n,r+1} = \prod_1^{p_{n,r+1}} T_j$  (each $T_j$ of length $1 - \frac{1}{p_{n,r+1}} \log \frac{\beta_r}{a_{n,r+1}}$ in $I_j$, $\forall j \in J_{n,r+1}$).  Thus
\[
	\mu(\bar{C}_{n,r+1}) = \left (1 - \frac{1}{p_{n,r+1}} \log \frac{\beta_r}{a_{n,r+1}} \right ) ^ {p_{n,r+1}}
\]
and $\mu(\bar{C}_{n,r+1}) \rightarrow \frac{a_{n,r+1}}{\beta_r}$ as $p_{n,r+1} \rightarrow \infty$.\\

\noindent We can thus take $p_{n,r+1}$ large enough so that $\frac{1}{2} \frac{a_{n,r+1}}{\beta_r} \leq \mu(\bar{C}_{n,r+1}) \leq \frac{a_{n,r+1}}{\beta_r}$ and that 
\[
	1 - \frac{1}{p_{n,r+1}} \log \frac{\beta_r}{a_{n,r+1}} > \frac{1}{2};
\]

\noindent If we take $A_{n,r+1} = \bar{B}_{n} \times \bar{E}_{n,r} \times \bar{C}_{n,r+1} \times I^{K_{n,r+1}^\prime}$, then  $\mu(A_{n,r+1}) \leq a_{n,r+1}$.\\

\noindent For each $j \in J_{n,r+1}$, again $S_j = \prod_{\substack{i \in J_{n,r+1}\\i \neq j}} T_i$.\\

\noindent For $J = K_{n,r} \cup (J_{n,r+1} - \{j\})$, we have $f_J(x) > \frac{1}{2}$ for all $x \in \bar{B}_{n} \times \bar{E}_{n,r} \times S_j \times I_j \times I^{K_{n,r+1}^\prime}$.\\

\noindent If $D_{n,r+1}$ is the union of the $p_{n,r+1}$ sets for each $j \in J_{n,r+1}$, it is clear that $\phi_{J_n}(x) > \frac{1}{2}$ in $D_{n,r+1}$.\\

\noindent However, we have $D_{n,r+1} = \bar{B}_{n} \times \bar{E}_{n,r} \times \bar{D}_{n,r+1} \times I^{K_{n,r+1}^\prime}$, where\\$\bar{D}_{n,r+1} = \bigcup_1^{p_{n,r+1}} (S_j \times I_j)$ and thus
\[
	\mu(\bar{D}_{n,r+1}) = \left (1 - \frac{1}{p_{n,r+1}} \log \frac{\beta_r}{a_{n,r+1}} \right ) ^ {p_{n,r+1}} + \left (\log \frac{\beta_r}{a_{n,r+1}} \right ) \left (1 - \frac{1}{p_{n,r+1}} \log \frac{\beta}{a_{n,r+1}} \right ) ^ {p_{n,r+1} - 1}.
\]

\noindent Taking into account the assumption that $\beta_r \geq \frac{1}{2}$, we can suppose $p_{n,r+1}$ is large enough that  $\frac{1}{2} \frac{a_{n,r+1}}{\beta_r} \log \frac{1}{a_{n,r+1}} \leq \mu(\bar{D}_{n,r+1}) \leq 2 \frac{a_{n,r+1}}{\beta_r} \log \frac{1}{a_{n,r+1}}$.\\

\noindent We deduce at once that $\frac{1}{2} a_{n,r+1} \log \frac{1}{a_{n,r+1}} \leq \mu(\bar{E}_{n,r} \times \bar{D}_{n,r+1}) \leq 2 a_{n,r+1} \log \frac{1}{a_{n,r+1}}$.\\

The recurrence on $r$ can thus continue just like this until we arrive at an $r$ such that $\delta_{n,1} + \delta_{n,2} + \dots + \delta_{n,r} > \frac{1}{2}$ and as $(a_n \log \frac{1}{a_n})$ diverges by assumption, there always exists a smaller $r$ having this property; it is this $r$ which we will take to be our $h_n$.

It is clear then that the $h_n$ sets $A_{n,1}, \dots, A_{n,h_n}$ are pairwise disjoint; so are the $h_n$ sets $D_{n,1}, \dots, D_{n,h_n}$.  Moreover, one has the union $D_n$ of $D_{n,r}$ has measure greater than $\frac{7}{16}$, and in all points of $D_n$ we have $\phi_{J_n}(x) > \frac{1}{2}$.  The union $A$ of all the sets $A_{n,r}$ thus answers the question.

\section{Consequence of Dieudonn\'e's example}

The preceding example allows us to respond in the negative to an analogous question concerning 
derivation bases\footnote{A net $\mathfrak{H}$ is a countable family of bounded, non-empty Borel sets such that:
\begin{itemize}
	\item $\{B: \mathfrak{H} \ni H \subset B\} \cap \mathfrak{H}$ is a finite family.
	\item	If $B_1, B_2 \in \mathfrak{H} \mbox{ and } B_1 \cap B_2 \neq \emptyset$, then either $B_1 \subset B_2 \mbox{ or } B_2 \subset B_1$.
\end{itemize}
Let $\mathfrak{M}$ be a $\sigma$-algebra on a set $X$, and let $E \subseteq X$ be fixed.  For each $x$ in $E$, let $(M_i(x))$ be a net.  The family of all $(M_i(x))$ forms a \emph{prebasis} $\mathcal{B}$.  Thus, the elements of $\mathcal{B}$ are converging sequences together with their convergence points.  We allow two or more points corresponding to the same sequence.  Let $\mathcal{D}$ be the family of all sets occurring in the sequences $(M_i(x))$ for all $x$ in $E$.  If we provide $\mathfrak{M}$ with a measure $\mu$ and if the sets of $\mathcal{D}$ are of finite measure, then $\mathcal{B}$ is a \emph{derivation basis}.}in the cube $\II^\infty$.  For the cube of finite dimension $\II^d$, a classical theorem of Vitali shows that the cubes $\Pi_n(x)$ with centre $x \in \II^d$ and sides $\frac{1}{n}$ form a 
derivation basis for the functions.  In particular, for all measurable, bounded $f$ defined in $\II^d$, the functions
\[
	g_n(x) = \frac{1}{\mu(\Pi_n(x))} \int_{\Pi_n(x)} f \, d\mu
\]
tend almost everywhere to $f(x)$ in $\II^d$.  One naturally wonders whether the following similar result is valid in $\II^\infty$: for all finite subsets $J$ of $\NN$ and all $n$, denote by $\Pi_{n,J}(x)$ the product of the cube with centre $x_J$ and with sides $\frac{1}{n}$ in $\II^J$ by $\II^{J^\prime}$, and for all measurable, bounded functions $f$ in $\II^\infty$, consider the function 
\[
	g_{n,J}(x) = \frac{1}{\mu(\Pi_{n,J}(x))} \int_{\Pi_{n,J}(x)} f \, d\mu;
\]
This function tends almost everywhere to $f(x)$ according to the ideal product $\NN \times F$ (with order relation $(n_1, J_1) \leq (n_2, J_2)$ signifies ``$n_1 \leq n_2$ and $J_1 \subset J_2$'').  The example constructed shows that the answer is no.  Indeed, let us write $\Pi_{n,J}(x) = \bar{\Pi}_{n,J}(x) \times \II^{J^\prime}$, where $\bar{\Pi}_{n,J}(x)$ is the cube with centre $x_J$ and side $\frac{1}{n}$, one has, with the notations above and by virtue of Fubini's theorem
\[
	\int_{\Pi_{n,J}(x)} f \, d\mu = \int_{\bar{\Pi}_{n,J}(x)} f_J \, d\mu;
\]
and $\mu(\Pi_{n,J}(x)) = \mu(\bar{\Pi}_{n,J}(x))$.  Recall that Vitali's theorem above shows that, for all $J$, one has 
\[
	\lim_{n \rightarrow \infty} g_{n,J}(x) = f_J(x) \mbox{ almost everywhere.}
\]
Since the sets of $F$ are countable, there thus exists in $\II^\infty$ a set of measure zero in the complement, of which one would have
\[
	\lim_{n \rightarrow \infty} g_{n,J}(x) = f_J(x) \mbox{ for all } J \in F.
\]
But then if $g_{n,J}(x)$ tended almost everywhere to $f(x)$ according to $\NN \times F$, the theorem of double limits \cite{NB} shows that $f_J(x)$ tends almost everywhere to $f(x)$ according to $F$, which is what we showed to be inaccurate.  Thus, the most general form fails, but maybe a form more suited to our situation will suffice.




\cleardoublepage

\chapter{Slowly Oscillating Functions}

\section{Jessen's Theorem}


Jessen showed that in countable cartesian product spaces, the infinite-dimensional integral is the limit (in the sense of $L^1$-norm) of the corresponding interval over the first $d$ spaces.

Jessen proved that as $n \rightarrow \infty$ and for almost every $x = (x_1, x_2, \ldots)$:
\begin{itemize}
	\item [(1)]	$\int \int \ldots f(x) \, dx_d dx_{d+1} \ldots \rightarrow f(x)$
	\item [(2)]	$\int \int \ldots \int f(x) \, dx_1 dx_2 \ldots dx_d \rightarrow \int_X f(x) dx$
\end{itemize}
where $X$ is the product of a countable sequence of measure spaces $X_1, X_2, \ldots, X_d, \ldots$, each of measure 1; and $f$ is a summable real-valued function of $X$.  As pointed out by examiner W. Jaworksi, Jessen's theorem is an immediate consequence of martingale theory (martingale convergence in the case of formula (1) and reverse martingale convergence in the case of formula (2)).  Martingales were not available to Jessen, but are a standard tool today.

Dorothy Maharam shows in her paper \cite{DM} how to extend (2) to the product of arbitrarily many coordinate spaces by taking the integrals over all finite subsets of the coordinate spaces.  As we saw in the previous chapter, Dieudonn\'e shows that such an extension is false for (1).  

The point is very subtle.  Jessen shows that $\int \int \ldots f_{[1,2,...d]}(x) \, dx_{d+1} \ldots \rightarrow f(x)$ not $\int \int \ldots f_J(x) \, dJ^{\prime} \ldots \rightarrow f(x)$ where $J \subset Fin$ and $Fin$ is the set of all finite subsets of $\NN$, and $J^{\prime}$ is its complement.  So the order of factors is important.  Maharam does show that an extension for (1) is possible if we use factors which are well-ordered and transfinite limits instead of directed limits.

Recall from the chapter \ref{C:Dieudonne} that Jessen works in $\II^\infty$ and for $x \in \II^\infty$, we can write $x = (x^\prime, x^{\prime \prime})$ where $x^\prime \in \II^J$ and $x^{\prime \prime} \in \II^{J^\prime}$.  Consider on each $\II_n$ (that is, each copy of $\II$), the Lebesgue measure, and denote by $\mu$ the product measure of the measurable sets of $\II^\infty$. We denote in the same way by $\mu_J$ the product measure on $\II^J$.  That being, say $f$ is a function defined on $\II^\infty$ and is integrable with respect to $\mu$.  According to the Fubini theorem, if $J$ is any subset of $\NN$ and $J^\prime$ is the complement, for almost all $x_J \in I^{J^\prime}$, the function $x_J \rightarrow f(x_J, x_{J^\prime})$ is integrable, the function
\[
	f_J(x) = \int_{I^{J^\prime}} f(x_J, x_{J^\prime}) d\mu_{J^\prime}
\]
defined almost everywhere in $I^J$ is integrable in this set, and one has
\[
	\int_{\II^\infty} f d\mu = \int_{I^J} f_J(x_J) d\mu_J.
\]

That being, the theorem by Jessen is as follows:
\begin{theo}
	If $(J_n)$ is an increasing sequence of finite subsets of $\NN$, the functions $f_{J_n}$ converge almost everywhere to f in $\II^\infty$ as $n$ goes to infinity.
\end{theo}

\section{Results on $\II^\infty$} \label{S:ResultsOnIinfty}

For continuous functions the Lebesgue density theorem is true in a very general sense.


A more general proof of theorem \ref{T:LebAtContinuousPt} is as follows:

\begin{theo} [\ddag] \label{T:LebMetricSpace}
	Let $(X, \rho, \mu)$ be a metric space with Borel probability measure and suppose $\mu$ has full support, that is, $\mbox{supp }(\mu) = X$.  Let $f:(X, \rho) \rightarrow \RR$ be a Borel function 	which is continuous at some point $x$.  Then 
	\[
		\displaystyle \lim_{\mu(V) \rightarrow 0} \frac{1}{\mu(V)} \int_V |f(x) - f(y)| d\mu(y) = 0
	\]
	in the following sense: for every $\epsilon > 0$, there is an open neighbourhood $V$ of $x$ such that whenever $B \subseteq 		V$ is a Borel subset of strictly positive measure,
	\[
		\frac{1}{\mu(B)} \int_{B} |f(x) - f(y)| d\mu(y) < \epsilon.
	\]
\end{theo}

\begin{proof}
	Let $x \in supp(\mu)$.  Since $f$ is continuous at $x$, given $\epsilon > 0$, $\exists V \ni x$ such that $\forall y \in V, 		|f(x) - f(y)| < \epsilon$.  This implies that for every Borel set $B$ with $\mu(B) > 0, \int_{B} |f(x) - f(y)| d\mu(y) < 					\epsilon \cdot \mu(B)$ and so $\frac{1}{\mu(B)} \int_{B} |f(x) - f(y)| d\mu(y) < \epsilon$ holds.
\end{proof}
\begin{rmk}
	Since $(X, \rho)$ is a metric space, we can always find a V as above.  Take as an example, $V_n = B(x, \frac{1}{n})$, the 			ball of radius $\frac{1}{n}$ around $x$.  These balls form a basis for the topology on X.  Taking $\mu \left (\bigcap_{n \in 		\NN} V_n \right ) = \lim_{n \rightarrow \infty} \mu(V_n) = 0$, we have a sequence of sets with their measures tending to 				zero.  Also, the proof of theorem \ref{T:LebMetricSpace} does not apply to $\RR^\infty$ since our measure $\lambda^\infty$ is 	not a probability measure.
\end{rmk}


Now the space $(\II^{\infty}, \rho)$ is a compact metric space when equipped with the metric which induces the product topology, for example,
\[
	\rho(x,y) = \sum_{n=1}^{\infty} 2^{-n} \frac{|x_n - y_n|}{1 + |x_n - y_n|}.
\] 
Let $\RR$ be equipped with the standard Euclidean norm.  Let $\pi_d: 	\II^{\infty} \rightarrow \II^{d}$ denote a projection, that is, if $x \in \II^{\infty}$ and $x = \{x_1,x_2,\ldots,x_{d-1},x_d,x_{d+1},\ldots\}$, then $\pi_d$ truncates $x$ to $x_d = \{x_1,x_2,\ldots,x_{d-1},x_d\}$.  Let $f: \II^{\infty} \rightarrow \RR$ be a continuous function with respect to the product topology, then by the Heine-Cantor theorem, $f$ is also uniformly continuous.  Thus, $\forall \epsilon > 0$, $\exists \delta > 0$ such that for all $x, y \in \II^{\infty}$ 	with $\rho(x,y) < \delta$, one has $|f(x) - f(y)| < \epsilon$.

\begin{theo} [\ddag] \label{T:UniformContinuity}
	Let $f:\II^\infty \to \RR$ be continuous.  For all $\epsilon > 0$, $\exists D$ such that $\forall d \geq D$, if $\pi_d(x) = \pi_d(y)$, then $|f(x) - f(y)| < \epsilon$.
\end{theo}

\begin{proof}
  Given any $\epsilon > 0$, choose a $D$ so that for $x, y 	\in \II^{\infty}$ one has:
	\[
		\sum_{n=1}^{D} 2^{-n} \frac{|x_n - y_n|}{1 + |x_n - y_n|} = 0;
		\sum_{n=D+1}^{\infty} 2^{-n} \frac{|x_n - y_n|}{1 + |x_n - y_n|} < \delta
	\]
	Now if $d \geq D$, and $\pi_d(x) = \pi_d(y)$, one has:
	\[
		\rho(x,y) = \sum_{n=1}^{D} 2^{-n} \frac{|x_n - y_n|}{1 + |x_n - y_n|} + \sum_{n=D+1}^{\infty} 2^{-n} \frac{|x_n - y_n|}{1 + 		|x_n - y_n|} \\
		< \delta
	\]
	Thus one has that $\forall \epsilon > 0$, $\exists D$ such that $\forall d \geq D$, if $\pi_d(x) = \pi_d(y)$, then \\$|f(x) - f(y)| < \epsilon$.
\end{proof}
In other words, what we have just proven is that, in the case where $f$ is continuous, the function $f_{d}$ which is obtained from $f$  by integration along the fibres (full spaces) in dimensions $d$ and greater, differs from $f$ by less than $\epsilon$.  For such $f$, the Lebesgue density theorem works.  Functions like these are what we want and need.  Let us denote by $S$ the above family of bounded functions $f$ on $\II^\infty$ such that for every $\epsilon > 0$ there exists a dimension $D$ such that for all $d \geq D$, if the truncations of two elements $x, y \in \II^\infty$ are equal then $|f(x) - f(y)| < \epsilon$.  These are functions which oscillate `slowly' in high dimensions along the fibres $\II^{J^\prime}$.  That is, as we go up to a certain high dimension, we can be sure that the function will not change by much.  
What kind of functions exist in $S$?  All continuous functions live in that space.  By the definition of the continuity for product topology we will be able to find a sufficiently high dimension that the function does not change much along the fibres. Also, given a function $g:\RR^d \rightarrow \RR$ with $g \in L^1(X, \mu)$, then $g \circ \pi_d \in S$.  This is a function which is exactly constant along the fibres in dimension $d$ and higher.  Not every $L^1(X, \mu)$ function satisfies this however.  
It can be possible that this is the only family of functions for which Lebesgue theorem holds.  The next natural question is: What is a natural metric for which this class is a complete metric space?

We shall use the same norm as in the space $L^{\infty}(X,\mu)$, that is, the \emph{essential supremum}.  The norm is defined as follows:
\[
	\|f\|_\infty = \inf\{C \geq 0: |f(x)| \leq C \mbox{ for almost every } x\}.
\]
Functions $f$ and $g$ are in the same equivalence class, $f \sim g$, if they are equal almost everywhere, and so belong to the same equivalence class.  Suppose $(X, \Sigma, \mu)$ is a space with measure, then for two functions $f,g: X \rightarrow \RR$ we have:
\[
	\esssup_{x \in X} f = \inf_{g: \mu(\{x:f(x) \neq g(x)\}) = 0} \sup_{x \in X} g(x)
\]
Let us try to equip $S$ with the above norm and see what happens.  
We end up with the following result:
\[
	f \in S \Longleftrightarrow f_d \stackrel{L^{\infty}}{\rightarrow} f.
\]
Put otherwise, the class $S$ consists of all functions $f$ which satisfy a stronger version of Jessen's theorem; the functions $f_d$ converge to $f$ not almost everywhere as with Jessen's theorem, but in $L^{\infty}$ norm, and this condition means that along the fibres they get smaller and smaller (uniformly smaller).  It means that every such function is measurable but not vice versa.

Define $\mbox{ess osc } f = \inf_{g \sim f }\sup |g(x) - g(y)|$.  A proof may emulate the following reasoning: Suppose $f \notin S$, then there exists an $\epsilon > 0$ such that for all $d$ there exists an $A \subseteq \II^d$ with $\mu(A) > 0$ and for all $x \in A$, $\mbox{ ess osc } (f(\pi_d^{-1}(x))) \geq \epsilon$.  If $||f - g||_{\infty} < \frac{\epsilon}{2}$ then for all $x \in A$, $\mbox{ ess osc }~(f(\pi_d^{-1}(x)))~\geq~\frac{\epsilon}{2}$.

This can be restated more accurately as:
Let $f:\II^{\infty} \rightarrow \RR$.  Then $\forall \epsilon, \exists d$ such that for almost every $x \in \II^d$ we have $\mbox{ess osc } (f(\pi_d^{-1}(x))) < \epsilon$.

\begin{theo} [\ddag]
	Let $S$ be the space of all bounded functions $f:~\II^{\infty}~\rightarrow~\RR$, with the property that $\forall 		\epsilon$, $\exists d$ such that $\forall d^{\prime} > d$,
\[
	\inf_{g \sim f} \sup_{\pi_d(x) = \pi_d(y)} |g(x) - g(y)| < \epsilon.
\]
Equip $S$ with the following norm: $||f|| = \mbox{ess sup } |f|$ where $\mbox{ess sup }_{x \in X} |f|$ is as defined earlier.  The space $S$ is complete.
\end{theo}

\begin{proof}
	Let $(f_n)$ be a Cauchy sequence in $S$ converging to some function $f$ almost everywhere.  It is enough to show that $f \in S$.  In other words, given a sequence $(f_n)$ satisfying $\forall \epsilon$, $\exists d$ such that $\forall d^{\prime} > d$, if $\inf_{g \sim f} \sup_{\pi_d(x) = \pi_d(y)} |g(x) - g(y)| < \epsilon$ and $(f_n)$ converges to $f$, then $f$ satisfies the above as well.  Let $\gamma > 0$ and choose for every $n$ a $g_n \sim f_n$ so that $\sup |g(x) - g(y)| < \epsilon + \gamma$.  Now, for almost every $x, y$ such that $\pi_d(x) = \pi_d(y)$ we have
	\[
	|f(x) - 	f(y)| < |f(x) - g_n(x)| + |g_n(x) - g_n(y)| + |g_n(y) - f(y)|
	\]
	and we already know that $f_n(x)$ converges to $f(x)$, thus $|f(x) - f_n(x)| < \epsilon$ almost everywhere and likewise $|g_n(y) - f(y)| < \epsilon$ almost everywhere.  Also, $|g_n(x) - g_n(y)| < \epsilon + \gamma$.  Thus $|f(x) - f(y)| < 3\epsilon + \gamma$.  Since $\gamma$ is arbitrary, we get $|f(x) - f(y)| \leq 3\epsilon$ and thus belongs to $S$ as well since $\epsilon$ is arbitrarily chosen.
\end{proof}

So far we have seen that examples of functions which are slowly oscillating include:
\begin{enumerate}
	\item continuous functions with regard to the product topology, and 
	\item	functions that factorise through projections.
\end{enumerate}

Now that we have positive results for $\II^\infty$, let us see if we can extend this to $\RR^\infty$ with the hope that they will survive.

\section{Results on $\RR^\infty$}

We begin this section by asking a simple, yet relevant question.  What does it mean when we say that $f: R^\infty \rightarrow R$ is integrable with regard to our infinite-dimensional Lebesgue measure?  Simply put, the integral exists and is finite, but for our purpose, let us examine it further.

\begin{itemize}
	\item	For bounded functions whose values are within some interval $-N$ to $N$ and whose domain of integration has finite measure, we subdivide this integral into small subintervals so that for every $i$, we consider a partition
				\[
					-N < -N + \frac{2N}{i} < -N + \frac{4N}{i} < \ldots < N.
				\]
				Then we form a Lebesgue integral sum,
				\[
					L_i(f) = \sum_{k=0}^{i-1} \lambda^\infty(f^{-1}(-N + \frac{2kN}{i},-N + \frac{2(k+1)N}{i}))\cdot(-N + \frac{2kN}{i})
				\]
				and we say that the bounded function is integrable if $\lim_{i \rightarrow \infty} L_i(f)$ exists and is finite.  That 					is
				\[
					\int f(x) d\lambda^\infty(x) = \lim_{i \rightarrow \infty} L_i(f) < \infty.
				\]
	
	\item If the domain of integration, $A$, is unbounded, we write it as the union of pairwise disjoint sets, $A_i$, each of finite measure and then 
	\[
		\int_A f \, d\lambda^\infty = \sum_i \int_{A_i} f \, d\lambda^\infty.
	\]
This sum does not depend on the choice of family as long as it is the union of disjoint sets of finite measure.
	
	\item	If the function is unbounded, then for every natural number $N$, we define a cut-off function $F_N$ as follows
				\[
					\forall N \in \RR, f_N := \min\{N, \max\{-N,f\}\}; \mbox{ and } \int f(x) d\lambda^\infty(x) = \lim_{N \rightarrow 							\infty} \int f_N(x) d\lambda^\infty(x)
				\]
				and we say that the function is integrable if this limit exists and is finite. 

\end{itemize}

Now, it is obvious that $\RR^{\infty}$ can be covered by an uncountably infinite family of cubes.  But $f$ being integrable means that the set of the points for which it is non-zero is contained in the union of countably many of these cubes.  For the simple reason that $\int f = \sum_\Pi \int_\Pi f|_{\Pi}$ and this sum cannot be uncountable.  Most of its terms should be zero for it to be sound.  So on most cubes the restriction of $f$ to them is zero.  Only on countably many cubes will the restriction be non-zero.  Thus, we can apply the theorem of this result to every cube and since there are countably many, it will survive.

Let us try to explain.

\begin{lem}	\label{L:BorelSetSigmaFinite1}
	Let $f:\RR^\infty \rightarrow \RR$ be Borel measurable and $f \geq 0$, and \\$\int_{\RR^\infty}~f(x)~d\lambda^\infty(x)~=~1$.  		Then $S = \{x: f(x) > 0\}$ has $\sigma$-finite measure. 
\end{lem}

\begin{proof}
	Let $S_n = \{x: f(x) \geq \frac{1}{n}\}, n = 1, 2, 3, \ldots, S = \cup_{n=1}^\infty S_n$.
	Let us look at $S_1$: 
	\[
		1 = \int_S f d\lambda^\infty \geq \int_{S_1} f d\lambda^\infty \geq \int_{S_1} 1 d\lambda^\infty = \mu(S_1)
	\]
	Similarly, the same argument works for $S_n$:
	\[
		1 = \int_S f d\lambda^\infty \geq \int_{S_n} f d\lambda^\infty \geq \int_{S_n} \frac{1}{n} d\lambda^\infty = \frac{1}{n} 				\mu(S_n)	
	\]
	This means that for every $n$, the measure of $S_n$ is finite.
\end{proof}
By the following lemma we can cover the set of all those points where $f(x) > 0$ by cubes and restrict the consideration of our function to the union of cubes.

\begin{lem}	[\ddag] \label{L:BorelSetSigmaFinite2}
	For a Borel set $A \subseteq \RR^\infty$, the following are equivalent:
	\begin{itemize}
		\item[(i)]	$A$ is $\sigma$-finite, that is, $A \subseteq \cup A_i$, where $\lambda^\infty(A_i) < \infty$ for all $i$ 
		\item[(ii)]	$A$ is contained in a union of a countable family of parallelepipeds of finite volume each.
	\end{itemize}
\end{lem}

\begin{proof}
	$(ii) \Rightarrow (i)$ is trivial.  $(i) \Rightarrow (ii)$ follows from definition of $\lambda^\infty$.  Let $\epsilon$ be any positive number, for example, $\epsilon = 1$.  Suppose $A \subseteq \cup_{i=1}^{\infty} A_i \mbox{ with } \lambda^\infty(A_i) < \infty$.  For each $i$, by definition there exists parallelepipeds $(C_{i,j})_{j=1}^{\infty}$ of finite measure such that $\lambda^\infty(A_i) \leq \sum_j \lambda^\infty(C_{i,j} \cap A_i) \leq \lambda^\infty(A_i) + \epsilon < \infty$.  Then $A \subseteq \cup_{i,j = 1}^{\infty} C_{i,j}$, each $C_{i,j}$ has finite volume.
\end{proof}
Note however that it may not always be possible to cover a Borel set $A \subseteq \RR^\infty$ by countably many cubes of side one, even in that case where the measure of $A$ is finite.  The following example shows this.

\begin{egg}
Let us suppose that we have a parallelepiped with sides of length greater than one which converge to one very fast so that the product exists and is finite.  Can it be covered with countably many cubes of finite volume?  No.  We can claim that by definition, a set of finite measure is contained in the union of countably many parallelepipeds of finite volume, but not necessarily with unit side.  Let $(a_n)~\downarrow~1 \mbox{ such that } \prod_{n=1}^\infty a_n < \infty$.  Let us examine the parallelepiped $[0, a_1] \times [0, a_2] \times \ldots$.  It requires uncountably many unit cubes to be covered.

Let $(C_1^n)_{n=1}^\infty$ be an infinite sequence of unit cubes.  We will define an $x \in [0, a_1] \times [0, a_2] \times \ldots$ for which $x \notin \cup_{n=1}^\infty C_1^n$.  For all $n, \mbox{ let } \mathcal{I}_j = [c_j, c_j+1]$ for some constant $c_j$ be the $j^{th}$ side of $C_1^n$.  Since $\mathcal{I}_j$ is a unit interval, looking at the $j^{th}$ interval of our parallelepiped, there exists an $x_j \in [0, a_j] \setminus [c_j, c_j+1]$.  This point is different from any $j^{th}$ coordinate of the $n^{th}$ cube so it misses all the cubes.  Thus we cannot cover the support of a measurable function with countably many cubes of side one.  We can do it with countably many parallelepipeds of finite volume though.
\end{egg}

\section{Peculiarities of $\RR^\infty$}

There are a few things to note about $\RR^\infty$ which are very peculiar.  The following theorem illustrates a problem that we came across.

\begin{theo} [$\ddag$]
	There are no functions $f:\RR^\infty \rightarrow \RR$ for which:
	\begin{itemize}
		\item[(i)]		$f \geq 0$
		\item[(ii)]		$\int_{\RR^\infty} f \, d\lambda^\infty = 1$
		\item[(iii)]	$f$ is continuous with regard to the product topology.
	\end{itemize}
\end{theo}

\begin{proof}
	Let $x_0 \in \RR^\infty$ with $f(x_0) > 0$.  Set $\epsilon = \frac{f(x_0)}{2} > 0$.  Since $f$ is continuous, it is continuous at $x_0$, and there exists a neighbourhood $V$ of $x_0$ in the product topology such that $\forall y \in V, |f(x_0) - f(y)| < \epsilon$.  Thus, for all $y$ in $V$, $f(y) \geq f(x_0) - \epsilon = \frac{f(x_0)}{2}$. By definition of the product topology, without loss of generality we can replace $V$ with a standard basic neighbourhood $V^\prime = \prod_{i=1}^d (a_i,b_i) \times \RR^{\{d+1, d+2, \ldots\}}$.  Now, notice that:
	\[
		1 = \int_{\RR^\infty} f(x) d\lambda^\infty(x) \geq \int_V f(x) d\lambda^\infty(x) \geq \int_V \epsilon d\lambda^\infty(x) = 		\epsilon \cdot \lambda^\infty(V) = +\infty
	\]
\end{proof}
\begin{rmk}
Moreover, the following stronger statement follows from the same argument.  If $f$ is a probability density function on $\RR^\infty$, then $f$ has to be discontinuous at every $x$ where $f(x) > 0$.  It can only be continuous at points where $f(x) = 0$.
\end{rmk}

\cleardoublepage

\chapter{Non-Density Theorem}

In this chapter we first discard an approach which once seemed plausible.  We then go on to give a simple example which explains the ideas of the main theorem.  Finally, we generalise the concepts of this example to prove  the main theorem of this work.

\section{Approach using Vitali systems}

Since the proof of the Lebesgue density theorem does not generalise to the infinite-dimensional case, what approach will we take?  Let us look at the most general form of the density theorem that is known as presented in \cite{SG} (section 10.3).  To state it, we need first to have the following definitions.

\begin{defn}
	A \emph{set function} is any function whose domain is a collection of sets and whose range is the (finite) real numbers.
\end{defn}

\begin{egg}
	We can take any measurable function $g$ and associate to it a set function $\phi$ as follows: $\phi(A) = \int_A g d\mu$ where $A$ belongs to the $\sigma$-algebra on which $\mu$ is defined.
\end{egg}

\begin{defn}
	A countably additive set function is a map $\phi: \mathfrak{B} \rightarrow \RR$ such that if $A_i, A_j \in \mathfrak{B} \mbox{ and } A_i \cap A_j = \emptyset \mbox{ for } i \neq j \mbox{ and } i,j \in \NN, $ then $\phi(\cup_{i=1}^\infty A_i) = \sum_{i=1}^\infty \phi(A_i)$, where $\mathfrak{B}$ denotes a $\sigma$-algebra.
\end{defn}

\begin{defn}
Let $X$ be a metric space equipped with a Borel measure $\mu$.  Suppose every set $\{x\}$ consisting of a single point is measurable with measure zero ($\mu(\{x\}) = 0$).  A \emph{Vitali system} for (X, $\mu$) is a family $\mathfrak{V}$ of Borel sets $E \subseteq X$ with
\begin{enumerate}
	\item Given any Borel set $E$, there are countably many sets $A_i, i = 1,\ldots,n,\ldots$ such that $E \subseteq 													\cup_{i=1}^{\infty}A_i \mbox{ and } \mu(\cup A_i \setminus E) < \epsilon$.
	\item	each $E \in \mathfrak{V}$ has a ``boundary'' $\partial E$ such that		
					\begin{enumerate}
						\item if $x \in E \setminus \partial E$ then all Vitali sets of sufficiently small measure containing $x$ are 														contained in $E \setminus \partial E$.
						\item	if $x \notin E \cup \partial E$ then all Vitali sets of sufficiently small measure containing $x$ are 															contained in $X \setminus (E \cup \partial E)$.
					\end{enumerate}
\end{enumerate}
\end{defn}

\begin{egg}
	The following are Vitali systems for the Euclidean space $\RR^d$ equipped with the Lebesgue measure:
	\begin{enumerate}
		\item The balls $B_x(\epsilon)$ of radius $\epsilon$ with centre $x \in \RR^d$ and $\epsilon > 0$.
		\item	All cubes.
	\end{enumerate}
\end{egg}

\begin{defn}
	Let $\phi(E)$ be a countably additive set function defined on a metric space $X$ (and hence on $\mathfrak{V}$).  Also, let $X$ be equipped with a Borel measure $\mu$.  Then, by the derivative of $\phi(E)$ at the point $x_0$ with respect to the Vitali system $\mathfrak{V}$ we mean the quantity
	\[
		D_{\mathfrak{V}} \phi(x_0) = \lim_{\epsilon \to 0} \frac{\phi(A_{x_0}(\epsilon))}{\mu(A_{x_0}(\epsilon))}
	\]
	(provided the limit exists), where $A_{x_0}(\epsilon)$ is any Vitali set of measure less than $\epsilon$ containing the point 	$x_0$.
\end{defn}
Put differently, differentiating $\phi$ at any point $x_0$ gives us $D_{\mathfrak{V}} \phi(x_0)$.

The following theorem is taken from \cite{SG} (section 10.3) and is the most general form of the theorem which we have come across.

\begin{theo} [Lebesgue-Vitali Theorem]
	Let $\mathfrak{V}$ be a Vitali system of Borel subsets of $X$ and let $\phi(E)$ be a countably additive set function on $X$.  	Then the derivative:
	\[
		D_{\mathfrak{V}} \phi(x_0) = \lim_{\epsilon \to 0} \frac{\phi(A_{x_0}(\epsilon))}{\mu(A_{x_0}(\epsilon))}
	\]
	exists almost everywhere.
\end{theo}

Let $\mathfrak{V}$ be a Vitali system of sets in $(X, \mu)$ where $\mu$ is $\sigma$-finite and $\sigma$-additive and let $f: X \rightarrow \RR$ be integrable.  Then this corollary follows immediately from the above theorem.

\begin{cor}
	Let $f$ be a measurable function on a metric space.  Then for almost every $x \in X$,
	\[
		\lim_{\epsilon \to 0} \frac{1}{\mu(A)} \int_A f(y) \, d\mu(y) = f(x)
	\]
	where $A \subseteq X$ is an element of a Vitali system containing $x$ and $\mu(A) < \epsilon$.
\end{cor}

Unfortunately, this approach cannot be made to work in $\RR^\infty$ for two main reasons.  First of all, the theorem uses a $\sigma$-finite measure, $\mu$, while the measure, $\lambda^\infty$, we are interested in is not $\sigma$-finite.  Second, there are no obvious candidates for Vitali systems.  Even in $\II^\infty$, the sets $\{[a_1, b_1] \times \ldots \times [a_n, b_n] \times \II^{J^\prime}\}$ do not form a Vitali system because the ``boundary'' condition does not work.  



A very relaxed approach to his first axiom is as follows: take a parallelepiped $C$ such that $0 < \lambda^\infty(C) < \infty, \forall x \in C, \exists \epsilon > 0 \mbox{ such that } \forall C^\prime \ni x, \lambda^\infty(C^\prime) < \epsilon \Rightarrow C^\prime \subseteq C$.  Take for example, $C = [0, 1]^\infty$.  If we take parallelepipeds, then we can always find a parallelepiped, $C^\prime$, of small measure which sticks out (thus violating the ``boundary'' condition) of our chosen parallelepiped of finite volume, $C$.  Notice that in  finite dimensions there is a very rigid dependence between the volume of a ball, or cube, and its radius, respectively the length of its side.  This is a one-to-one correspondence and by changing the volume we can make the radius as small as we wish.

In infinite dimensions, there is only one cube of finite dimension.  The cubes in $\RR^\infty$ have the property that their volumes are either 0, 1, or +$\infty$.  This is why we are forced to use parallelepipeds.  But parallelepipeds do not form a Vitali system in $\RR^\infty$ as we just saw.

Let us thus abandon the Vitali approach and recall that parallelepipeds do work in $\RR^d$.  In \cite{HP} (Part I, Chapter V) a theorem by Jessen, Marcinkiewicz and Zygmund shows that the density theorem holds in $\RR^d$ using a parallelepiped basis.  Their theorem is stated as follows:

\begin{theo}
The interval basis $\mathfrak{I} = [\mathcal{I}, \delta]$ in $\RR^d$ derives the Lebesgue integral of each measurable function $f$ for which the function $|f|(\log^+ |f|)^{m-1}$ is Lebesgue integrable over the open cube
\[
	Q_0 = \{x: 0 < x_i < 1, i = 1,2,\ldots,d\},
\]
and the $\mathfrak{I}$-derivative coincides with $f$ except on a set of Lebesgue measure zero.
\end{theo}

Here, $\mathcal{I}$ is the family of closed non-degenerate $d$-dimensional parallelepipeds
\[
I = \{x: \alpha_i \leq x_i \leq \beta_i, i = 1,2,\ldots,d\}
\]
and $\alpha_i < \beta_i \mbox{ for } i = 1,2,\ldots,d$.  We will not analyse it but the condition on $f$ is satisfied by all bounded functions for example.  In different words, the above theorem states that the derivative of any Lebesgue integral of a measurable function $f$, satisfiying certain conditions, can be calculated using the interval basis and the derivative will coincide with $f$ except on a null set.

\begin{rmk}
Interestingly enough, the same source (\cite{HP}, p. 104) shows that in the more general case, where the basis consists of 			all rectangular parallelepipeds whose sides may or may not be parallel to the coordinate axes, the density theorem does not 	even hold in $\RR^2$.
\end{rmk}



\section{Interpretation of Density}

It is important to realize that since we are dealing with parallelepipeds, there are two ways in which we can interpret density.  The reason for this is that we can be given a parallelepiped $\Pi$ with centre $x$ and measure $0 < \lambda^\infty(\Pi) < \infty$ and another parallelepiped $\Pi^\prime$ with centre $x$ and measure $\lambda^\infty(\Pi^\prime) < \lambda^\infty(\Pi)$ then one of two possible situations holds:

\begin{itemize}
	\item	The first situation is the one which comes to mind immediately and that is $\Pi^\prime \subseteq \Pi$.  For this case, the definition of a Lebesgue point is as follows:
	\begin{defn}
Call $x \in \RR^\infty$ a Lebesgue point for $f:\RR^\infty \rightarrow \RR$ if $\forall \epsilon > 0, \exists \Pi$, a parallelepiped centred at $x$, $\lambda^\infty(\Pi) > 0$, such that $\forall \Pi^\prime$ centred at $x$, $\lambda^\infty(\Pi^\prime) > 0, \Pi^\prime \subseteq \Pi$,
	\[
		|\frac{1}{\lambda^\infty(\Pi^\prime)} \int_{\Pi^\prime} f(y) d\lambda^\infty(y) - f(x)| < \epsilon.
	\]
\end{defn}
	\item	The second situation stems from the fact that the geometry of a parallelepiped is not determined by its measure and so even though $\lambda^\infty(\Pi^\prime) < \lambda^\infty(\Pi)$, it does not mean that $\Pi^\prime$ is contained in $\Pi$.  There may be dimensions in which $\Pi^\prime$ sticks out of $\Pi$.  As a simple example in $\RR^2$, take the cube $Q$ of side 2 units.  It's area is 4 units$^2$ of course.  Now, take a rectangle $R$ with the same centre as $Q$, but whose length is 3 and width is 1.  It's area is 3 units$^2$ but cannot be contained in $Q$.  In this case, the definition of a Lebesgue point has to be modified to the following stronger statement:
\begin{defn}
Call $x \in \RR^\infty$ a Lebesgue point for $f:\RR^\infty \rightarrow \RR$ if $\forall \epsilon > 0, \exists 									\delta$ such that for every parallelepiped centred at $x$ with $\lambda^\infty(\Pi) < \delta$ 
\[
	|\frac{1}{\lambda^\infty(\Pi)} \int_{\Pi} f(x) d\lambda^\infty(x) - f(x)| < \epsilon.
\]
\end{defn}
\end{itemize}

The second condition is a stronger one as if the density exists in that case, then it also exists in the first case and furthermore, they are equal.  However, it is not known whether the density exists for the second case even in finite dimensions and for this reason we will only work with the first case where the restriction is by geometry also and not only volume.


\section{Example: Calculating the density of the cube}

We need to realise some facts about density.  First, they do not always exist for measurable functions at all points, even in $\RR^d$.  The following example illustrates this.

\begin{egg}
The function
\(
	f(x) = \left \{
		\begin{array}{cl}
			1	&	\mbox{ if } x \in (2^{-(n + 1)}, 2^{-n}] \mbox{ for $n = 0, 2, 4, \ldots$ } \\
			0	&	\mbox{ otherwise }
		\end{array}
	\right.
\) has integral $\sum_{n=1}^{\infty} 2^{-(2n-1)} = \frac{2}{3}$ over the interval $[0,1]$, and thus integral $\frac{1}{3}$ over the interval $[-1,1]$.  What is the density of $f$ at $x = 0$?
\[
	\lim_{n \to \infty} \frac{1}{|[-2^{-n}, 2^{-n}]|} \int_{[-2^{-n}, 2^{-n}]} f(y) dx = ?
\]
Realise that if we take a ball of decreasing radius about $x = 0$, the density will never exist as the above limit does not exist.  The limit does not converge because it keeps fluctuating and does so even more as the ball gets smaller.
\end{egg}

Can we state a reasonable analogue of the Lebesgue density theorem for $\lambda^\infty$?  Let us begin with an example which contains the ideas of the main theorem.

\begin{egg}
Let us compute the Lebesgue density of $A = \II^\infty \subseteq \RR^\infty$.
\begin{align}
	density(A,x)	&= \lim_{\substack{x \text{ is the centre of } \Pi\\0 < \lambda^\infty(\Pi) < \epsilon\\\epsilon \rightarrow 0}} \frac{1}{\lambda^\infty(\Pi)} \int_\Pi \chi_A d\lambda^\infty \notag \\
								&= \lim_{\epsilon \to 0} \frac{\lambda^\infty(\Pi \cap A)}{\lambda^\infty(\Pi)} \notag
\end{align}

First, let us study the case where $x \in \II^\infty$.  Consider the set of sequences $S = \{x \in \II^\infty: \exists D, \forall d > D, x_d \geq \frac{1}{4}\}$.  This set has $\lambda^\infty$-measure zero.  Indeed, let $S = \cup_D S_D \mbox{ where } S_D = \{x: \forall d > D, x_d \geq \frac{1}{4}\} \subseteq [0,1] \times \ldots \times [0,1] \times [\frac{1}{4}, 1] \times [\frac{1}{4}, 1] \times \ldots$, which has measure zero.  We draw the conclusion by using the $\sigma$-additivity of $\lambda^\infty$ that for almost every $x \in \II^\infty, \forall D, \exists d > D, x_d < \frac{1}{4}$.

Now, let $x \in \II^\infty \setminus S$.  Take $C_1$ to be the unit cube centred at $x$.  Clearly, one has 
\[
C_1 \cap \II^\infty \subseteq [0,1] \times \ldots \times  [0, \frac{3}{4}] \times \ldots \times [0, \frac{3}{4}] \times \ldots
\]
where the interval $[0, \frac{3}{4}]$ occurs and infinite number of times.  This is a set that has measure zero.  It follows that 
\(
density(\II^\infty, x)  = 0.
\)

Finally, consider the case where $x \notin \II^\infty$.  Then there is a coordinate, $x_i$, which sits outside $\II^\infty$.  
Take a parallelepiped centred at $x_i$ which is so small in the $i^{th}$ dimension that it does not touch $\II^\infty$.  Since this parallelepiped is disjoint from $\II^\infty$, the measure of their intersection is zero.  
Overall, we conclude that
\[
\frac{\lambda^\infty(\Pi \cap \II^\infty)}{\lambda^\infty(\Pi)} = \frac{0}{1} = 0
\]
and the density of the cube $\II^\infty$ is zero at $\lambda^\infty$-almost every point of $\RR^\infty$.
\end{egg}

\section{Main Theorem}

Before we extend the previous ideas to a more general situation, let us introduce one more concept that is needed to prove our main result.  

\begin{defn} [$\ddag$]
Let $\Pi$ be a parallelepiped with centre $x$.  That is, let 
\(
\Pi~=~\prod_{i=1}^\infty [a_i, b_i]
\)
and
\(
	x = \mbox{centre}(\Pi) = \left ( \frac{a_i + b_i}{2} \right )_{i=1}^\infty.
\)
Let $0 \leq \delta \leq 1$.  The \emph{$\delta$-core} of $\Pi$ is the set:
\[
	\mbox{core}_\delta(\Pi) =\bigcup_{D=1}^\infty \left ( \prod_{i=1}^D [a_i, b_i] \times \prod_{i=D+1}^\infty [\frac{a_i + b_i}{2} - \frac{\delta \cdot (b_i - a_i)}{2}, \frac{a_i + b_i}{2} + \frac{\delta \cdot (b_i - a_i)}{2}] \right ).
\]
\end{defn}

The core consists of all sequences whose coordinates eventually end up close to the centre of the parallelepiped, and stay there forever.   Note that if $\delta = 1$, we get back to our ``parent'' parallelepiped.  


\begin{lem} [$\ddag$] \label{L:NullCore}
$\lambda^\infty(\mbox{core}_\delta(\Pi)) = 0$ if $\delta < 1$, provided $\lambda^\infty(\Pi) < \infty$.
\end{lem}

\begin{proof}
It is enough to prove that for each $D$,
\[
\prod_{i=1}^D [a_i, b_i] \times \prod_{i=D+1}^\infty [\frac{a_i + b_i}{2} - \frac{\delta \cdot (b_i - a_i)}{2}, \frac{a_i + b_i}{2} + \frac{\delta \cdot (b_i - a_i)}{2}]
\]
has measure zero because the core is the union of countable many sets and if each set has measure zero, then by the $\sigma$-subadditivity of measure, their union is also zero.

So:
\begin{align}
	\lambda^\infty & \left ( \prod_{i=1}^D [a_i, b_i] \times \prod_{i=D+1}^\infty [\frac{a_i + b_i}{2} - \frac{\delta \cdot (b_i - a_i)}{2}, \frac{a_i + b_i}{2} + \frac{\delta \cdot (b_i - a_i)}{2}] \right ) \notag \\
	&= \prod_{i=1}^D (b_i - a_i) \times \prod_{i=D+1}^\infty \delta \cdot (b_i - a_i) \notag \\
	&= \prod_{i=1}^D l_i \times \prod_{i=D+1}^\infty \delta \cdot l_i \notag
\end{align}
But $(l_i)_{i \to \infty} \to 1$ (see lemma \ref{L:SidesApproachOne}), so as $i \to \infty$ the lengths of the sides of the core converge to $\delta$.  But $\delta < 1$ and so $\prod_{i=D+1}^\infty \delta = 0$.
\end{proof}
We now arrive at the main theorem of this work.

\begin{theo} [$\ddag$] \label{T:NonDensity}
Let $f:\RR^\infty \rightarrow \RR$ be a measurable function such that $\int_{\RR^\infty}{f \, d\lambda^\infty}~=~1$ and $f \geq 0$.  Then for almost every $x$,
\[
\lim_{\substack{0 < \lambda^\infty(\Pi) < \infty\\\lambda^\infty(\Pi) \rightarrow 0}} \frac{1}{\lambda^\infty(\Pi)} \int_{\Pi} f(y) d\lambda^\infty(y) = 0.
\]
\end{theo}

\noindent [Proof Idea] Since $\int_{\RR^\infty} f \, d\lambda^\infty = 1$, the set of points where $f$ is non-zero, call it $A$, is $\sigma$-finite and so $A \subseteq \cup_{i=1}^\infty \Pi_i$ where $\Pi_i$ are parallelepipeds with $\lambda^\infty(\Pi_i) < \infty$. 
For each parallelepiped, the measure of the points at the $\delta$-core (for $\delta < 1$) is zero.  Since a countable union of negligible sets is still negligible, the measure of the union of the cores of all our parallelepipeds is still zero. Now, take a point which is not at the union of the cores and take the unit cube, $C_1$, around it.  Then $\lambda^\infty(C_1 \cap \Pi_i) = 0$ and 
%
so,
\[
	\frac{1}{\lambda^\infty(C_1)} \int_{\{y \in C_1;f(y) > 0\}} f(y) \, d\lambda^\infty(y) = 0
\]
and $A$ has density zero almost everywhere.

\begin{proof} 
Let $A = \{x: f(x) > 0\}$, then $A \subseteq \cup_{j=1}^\infty \Pi_j$ (see lemmas \ref{L:BorelSetSigmaFinite1} and \ref{L:BorelSetSigmaFinite2}).  Set $\delta = \frac{1}{2}$.  We have that $\lambda^\infty(\cup_{j=1}^\infty \mbox{ core}_\frac{1}{2} (\Pi_j)) = 0$ (see lemma \ref{L:NullCore}).  Without loss of generality, suppose $x \in \RR^\infty \setminus \{\cup_{j=1}^\infty \text{ core}_\frac{1}{2} (\Pi_j)\}$.  Denote by $C_1$ the unit cube centred at $x$.  Fix an arbitrary $j$.  We know that $x \notin \mbox{core}_\frac{1}{2} (\Pi_j)$ and this means that
\[
\left \vert \left \{ i: | x_i - \frac{b_i + a_i}{2} | > \frac{l_i}{4} \right \} \right \vert = \infty,
\]
where $l_i$ is the length of the sides of the parallelepiped in the $i^{th}$ dimension.  Set $\epsilon = \frac{1}{8}$ and choose $D$ so that $\forall i > D$, $\frac{7}{8} < l_i < \frac{9}{8}$.  In particular, for an infinite set of $i$'s we have both $\frac{7}{8} < l_i < \frac{9}{8}$ and $| x_i - \frac{b_i - a_i}{2} | > \frac{l_i}{4}$ being satisfied.

Let us deduce that $\lambda^\infty(C_1 \cap \Pi_j) = 0$.  We do so by bounding the length of the intersection of the sides of $C_1$ and $\Pi_j$.  The worst cases are:
\begin{enumerate}
	\item when $x$ comes close to the centre of the $i^{th}$ interval and
	\item	when the length of the interval is at its longest
\end{enumerate}
because in either case the intersection can be large.  So take the longest possible length which is $l_i = \frac{9}{8}$, and the $x$ closest to the centre.  This means that $x$ would be at a distance of $\frac{1}{4} \times \frac{9}{8} = \frac{9}{32}$ which is midway between the start and centre of the interval.  Now, let us calculate the overlap with the unit interval whose centre would be at that point $x$.  We see that the overlap has length $\frac{1}{2} + \frac{9}{32} = \frac{25}{32} < 1$.  Recall that this happens infinitely many times.

So
\[
\lambda^\infty(C_1 \cap \Pi_j) \leq 1 \times \ldots \times \frac{25}{32} \times \ldots \times \frac{25}{32} \times \ldots = 0.
\]
And as this happens for all $j$, it follows that
\[
\lambda^\infty(C_1 \cap \cup_{j=1}^\infty \Pi_j) \leq \sum_{j=1}^\infty \lambda^\infty(C_1 \cap \Pi_j) = 0
\]
and also,
\[
	\frac{1}{\lambda^\infty(C_1)} \int_{C_1} f(y) \, d\lambda^\infty(y) = 0.
\]
\end{proof}
The next corollary follows directly from the preceding theorem.

\begin{cor} [$\ddag$]
Let $A \subseteq \RR^\infty$ be a Borel set with $0 < \lambda^\infty(A) < \infty$.  Then for almost every $x$, $density(A, \lambda^\infty) = 0$, that is,
\[
	\lim_{\substack{\Pi \text{ with centre } x\\0 < \lambda^\infty(\Pi) < \infty\\\lambda^\infty(\Pi) \to 0}} \frac{\lambda^\infty(\Pi \cap A)}{\lambda^\infty(\Pi)} = 0.
\]
\end{cor}



At this point, we have come to the end of this work.  Essentially, we have shown that the approach used is not plausible.  It is not to be discarded though as there have been positive results.  However, it does not mean that positive results regarding densities may not be obtained by using a different measure.

\cleardoublepage





\newpage

\end{document}